\documentclass[11pt]{article}
\usepackage[utf8]{inputenc}
\usepackage{amsthm,amsmath,amssymb,amsfonts}
\usepackage{color}
\usepackage{a4wide}
\newcommand{\pder}[2]{\frac{\partial #1}{\partial #2}}

\newcommand{\dx}{\,{\rm d}x}

\newcommand{\Div}{{\rm div}\,}
\newtheorem{thm}{Theorem}
\newtheorem{lem}[thm]{Lemma}

\newtheorem{df}{Definition}
\newtheorem{rmk}{Remark}

\newcommand{\vr}{\varrho}

\newcommand{\vrn}{\vr_n}
\newcommand{\vtn}{\vt_n}

\newcommand{\vun}{\vu_n}

\newcommand{\vt}{\vartheta}

\newcommand{\vu}{\vc{u}}

\newcommand{\vc}[1]{{\bf #1}}
\newcommand{\vcg}[1]{{\pmb #1}}
\newcommand{\F}[1]{$\mathbb{#1}$}

\newcommand{\tn}[1]{\mbox {\F #1}}

\newcommand{\dS}{\,{\rm d}S}

\newcommand{\intO}[1]{\int_{\Omega} #1\dx}

\newcommand{\R}{\mathbb{R}}
\newcommand{\Rr}{\mbox{\FFF R}}
\font\FFF=msbm10 scaled 700
\newcommand{\N}{\mathbb{N}}

\begin{document}

\title{Steady compressible Navier--Stokes--Fourier system with general temperature dependent viscosities I: density estimates based on Bogovskii operator}
\author{Ond\v{r}ej Kreml $^1$, Tomasz Piasecki$^2$, Milan Pokorn\' y$^3$ and Emil Sk\v{r}\'\i\v{s}ovsk\'y$^3$}
\maketitle

\bigskip

\centerline{ $^{1}$ Institute of Mathematics, Czech Academy of Sciences}
\centerline{\v{Z}itn\'a 25, 115 67 Praha 1, Czech Republic}
\centerline{e-mail: {\tt kreml@math.cas.cz}}

\centerline{ $^{2}$ Institute of Applied Mathematics and Mechanics, University of Warsaw}
\centerline{ul. Banacha 2, 02-097 Warszawa, Poland}
\centerline{e-mail: {\tt tpiasecki@mimuw.edu.pl.}}

\centerline{$^{3}$ Charles University, Faculty of Mathematics and Physics}
\centerline{Mathematical Inst. of Charles University, Sokolovsk\' a 83, 186 75 Prague 8, Czech Republic}
\centerline{e-mail: {\tt pokorny@karlin.mff.cuni.cz}, {\tt emil.skrisovsky@gamil.com}}
\vskip0.25cm


\smallskip


\begin{abstract}

The aim of this paper is to reconsider the existence theory for steady compressible Navier--Stokes--Fourier system assuming more general condition of the dependence of the viscosities on the temperature in the form $\mu(\vt)$, $\xi(\vt) \sim (1+\vt)^\alpha$ for $0\leq \alpha \leq 1$. This extends the known theory for $\alpha=1$ from \cite{NP1} and improves significantly the results for $\alpha =0$ from \cite{MP2}. This paper is the first of a series of two papers dealing with this problem and is connected with the Bogovskii-type estimates of the sequence of densities. This leads, among others, to the limitation $\gamma >\frac 32$ for the pressure law $p(\vr,\vt) \sim \vr^\gamma + \vr\vt$. The paper considers both the heat-flux (Robin) and Dirichlet boundary conditions for the temperature as well as both the homogeneous Dirichlet and zero inflow/outflow Navier boundary conditions for the velocity. Further extension for $\gamma >1$ only is based on different type of pressure estimates and will be the content of the subsequent paper.  

\end{abstract}

\section{Introduction}
\label{1}

We consider the steady variant of the compressible Navier--Stokes--Fourier (NSF) system in a bounded domain $\Omega \subset \R^3$ which in the case of smooth solutions read as
\begin{equation} \label{NSF_1}
\begin{aligned}
\Div(\vr\vu) &=0 \\
\Div(\vr\vu\otimes\vu) -\Div \tn{S} + \nabla p &= \vr \vc{f} \\ 
\Div(\vr\vu E) +\Div\vc{q} -\Div (\tn{S} \vu) + \Div(p\vu) &= \vr \vc{f}\cdot \vu.
\end{aligned}
\end{equation} 
Above, the first equation represents the balance of mass (and it is usually called the equation of continuity), the second one represents the balance of linear momentum and the last one the balance of total energy. On the level of smooth solutions, the last equality (together with the first two equations) is equivalent with the balance of internal energy which is not suitable in the case of weak solutions (however, see \cite {BdV} in the case of regular solutions or even \cite{MP1} in a very particular case of more regular weak solutions) or with the entropy equality discussed below; the entropy will play an important role in our analysis.

Above, $\vr$: $\Omega \to \R^+_0$ stays for the density of the fluid, $\vu$: $\Omega \to \R^3$ for the velocity field. The third unknown function is the temperature $\vt$: $\Omega \to \R^+$. The tensor valued function (symmetric) $\tn{S}= \tn{S}(\vt,\nabla \vu)$ is the stress tensor, the scalar function $p=p(\vr,\vt)$ is the pressure of the fluid, the vector $\vc{q} = \vc{q}(\vt,\nabla \vt)$ is the heat flux and the scalar $E = e+\frac 12 |\vu|^2$ is the specific total energy. The specific internal energy $e=e(\vr,\vt)$ as well the stress tensor, the pressure and the heat flux are given functions of the unknowns and their derivatives. Finally, the given vector $\vc{f}$: $\Omega \to \R^3$ plays the role of the volume force.

System \eqref{NSF_1} describes the steady flow of heat-conducting compressible fluid. We need to specify the boundary conditions on $\partial \Omega$. In this paper we will consider for the velocity either the homogeneous Dirichlet conditions 
\begin{equation} \label{vel_Dir}
\vu = \vc{0}  
\end{equation}
or the Navier (partial or full slip) boundary conditions
\begin{equation} \label{vel_Nav}
\vu\cdot \vc{n} = 0, \qquad (\tn{S}\vc{n})\cdot \vcg{\tau} + \sigma \vu\cdot \vcg{\tau} = 0,
\end{equation}
where $\vc{n}$ is the normal vector to $\partial \Omega$ and $\vcg{\tau}$ is the tangent vector (there are two linearly independent tangent vectors at each point of $\partial \Omega$ in our situation) and the friction coefficient $\sigma \geq 0$. For the temperature we either assume the heat flux to be proportional to the difference of the temperature $\vt$ inside and $\vt_0$ outside, i.e.  
\begin{equation} \label{temp_heatflux}
\vc{q} \cdot \vc{n} = L(\vt-\vt_0)
\end{equation}
for $L > 0$ being a given constant, or we assume the temperature to be prescribed on the boundary $\partial \Omega$ (inhomogeneous Dirichlet boundary condition)
\begin{equation} \label{temp_Dir}
\vt = \vt_D.
\end{equation}
Finally, have to specify the total mass of the fluid
\begin{equation} \label{tot_mass}
\int_\Omega \vr \dx = M >0.
\end{equation}

While the existence of strong solutions (for small data or data close to some special solution) has been known since eighties of the last century (see \cite{BdV}) and the technique was further developed including inflow/outflow conditions in, e.g., \cite{KK,MPi,PP1} or \cite{Zh}, the situation with solutions for large data appeared to be more complex. A certain intermediate result can be found in \cite{AMP}, but it is still rather a small data result. The seminal monograph \cite{Li2} contains some ideas for this problem, however, under the assumption that not the total mass, but a certain $L^p$ norm of the density is given with $p$ sufficiently large. The first treatment of the problem stated above, with the Navier boundary condition for the velocity and the heat flux boundary condition for the temperature was performed in \cite{MP1}; this paper is based on ideas developed in \cite{MP_01} and \cite{MP_02} for the steady compressible Navier--Stokes equations.  The authors were able to show existence of solutions with bounded density and almost bounded velocity gradient which also fulfilled the weak formulation of the internal energy balance. The price to pay was the necessity to have relatively large $\gamma$ in the pressure law. Leaving the nice properties of the density and velocity it was possible to extend the range of pressure laws as well as to include the homogeneous Dirichlet boundary conditions for the velocity in \cite{MP2}. However, the a priori estimates for this problem were still based on the linear momentum balance combined with the entropy inequality which on one hand leads to $L^2$-estimates of the velocity gradient, on the other hand to a relatively high value of $\gamma$. 

A new idea appeared in  papers \cite{NP1}, \cite{NP2} (and \cite{NP3} for the two-dimensional case), where the main estimates were read from the entropy inequality combined with integrated version of the total energy balance. The price to pay here was the assumption that the viscosities grow linearly with respect to the temperature since the estimates are in this case optimal. A certain attempt to relax this assumption to more general dependences was  paper \cite{KNP}, however, the terms coming from the radiation effect included in the paper lead to suboptimal results with respect to the compressible NSF system itself. Note also that there are two possibilities how to estimate the density: either solely by the application of Bogovskii operator (applied in \cite{NP1}) which leads to a lower bound on the parameter ($\gamma>\frac 32$ for linear growth in the viscosity with respect to the temperature, this approach was also used in \cite{KNP}) or combining the Bogovskii operator estimates with estimates of the type  $\sup_{x_0\in\Omega} \int_\Omega \frac{p(\vr,\vt)}{|x-x_0|^\beta}\dx$ for some $1\geq \beta >0$. This method, based on ideas developed by Plotnikov and Sokolowski (see \cite {PS1}) and extended by Frehse, Steinhauer and Weigant (see \cite{FSW}) was applied for the NSF system in  \cite{JNP} or \cite{NP2}. The main idea here is to use test functions of the type $\frac{x-x_0}{|x-x_0|}$ in the momentum equation; if $x_0$ is far from the boundary, this function is multiplied by a suitable cut-off function, if it is close to or on the boundary $\partial \Omega$, certain modifications are needed. Some further results and ideas can be found in the overview paper \cite{MPZ}. Finally, let us mention the Dirichlet boundary condition for the temperature. Based on the idea of the ballistic energy inequality applied in the evolutionary case by Feireisl and Chaudhuri \cite{FC}, the steady case was considered in \cite{P}. 

The aim of this paper is to apply the former in the case of all above mentioned boundary conditions for sublinear growth of viscosities with respect to the temperature, including the case of constant viscosities. The results in the case of homogeneous Dirichlet and the Navier boundary condition do not differ at all (which will not be the case of the second paper, where the latter will be applied), however, they differ from the situation above for the Dirichlet boundary condition for the temperature due to the use of the ballistic energy inequality instead of the the entropy one. This paper deals with the situation when solely the Bogovskii operator estimated for the density are used. The other situation based on the special test functions presented above will be the content of a subsequent paper. Hence, here, we consider the case when $\gamma$ is not too close to 1 (in fact, it will be always larger than $\frac 32$) while the other approach allows to deal with $\gamma$ close to $1$.     

The paper is structured as follows. The next section contains the precise statement of assumptions on the constitutive assumptions for the model introduced above as well as  the definitions of weak solutions in all possible situations and the statement of the results of the paper. Section \ref{3} recalls several auxiliary results needed later. The last two sections include the proofs of the main result; first we consider the heat flux boundary conditions and finally the Dirichlet boundary condition. 

In what follows, a standard notation for the Lebesgue and Sobolev spaces will be used ($L^p(\Omega)$ or $W^{k,p}(\Omega)$, respectively) endowed with the standard norms denoted by $\|\cdot\|_p$ and $\|\cdot\|_{k,p}$, respectively. The generic constant will be denoted by $C$ and its value may change in the same formula or even on the same line.  

\section{Constitutive relations. Main results}
\label{2}

In this section, we first introduce the constitutive relations we shall deal with below, then recall the definitions of our solutions and finally formulate the main results of our paper.

\subsection{Constitutive relations}
\label{2.1}

We take as unknown functions the density $\vr$, the temperature $\vt$ and the velocity $\vu$. All other thermodynamic and mechanic quantities introduced below will be given functions of these unknown functions. 

First, we assume the fluid to be compressible Newtonian which results into
\begin{equation} \label{stress}
\tn{S} = \tn{S}(\vt,\nabla \vu) = \mu(\vt) (\nabla \vu + \nabla \vu^T-\frac 23 \Div \vu \tn{I}) + \xi (\vt) \Div \vu \tn{I} 
\end{equation}
where the shear viscosity $\mu(\vt) >0$ is a Lipschitz continuous function, the bulk viscosity $\xi(\vt) \geq 0$ is a continuous function and there exist constants $C_1$ and $C_2$ such that
\begin{equation} \label{visc}
\begin{aligned}
C_1 (1+\vt)^\alpha &\leq \mu(\vt) \leq C_2(1+\vt)^\alpha \\
0 &\leq \xi(\vt) \leq C_2(1+\vt)^\alpha 
\end{aligned}
\end{equation}
for some $\alpha \in [0,1]$. The fact that we also consider $0 \leq \alpha <1$ is the main novelty of the paper with respect to \cite{NP1}, \cite{NP2} or \cite{P}. On the other hand, papers \cite{MP1} and \cite{MP2} worked with constant viscosities and we will comment later on the improvements of the results from this paper with respect to these older papers. 

The heat flux is assumed in the form of the Fourier law
\begin{equation} \label{heat_flow}
\vc{q}=\vc{q}(\vt,\nabla \vt) = -\kappa(\vt)\nabla \vt,
\end{equation}
where the heat conductivity $\kappa(\vt)>0$ is a continuous function and there are positive constants $C_3$ and $C_4$ such that 
\begin{equation} \label{heat_cond}
C_3 (1+\vt)^m \leq \kappa (\vt) \leq C_4 (1+\vt)^m
\end{equation}
for a parameter $m>0$.

The last parameter which plays a role in our analysis is connected with the pressure. We assume a generalization of the pressure law for the monoatomic gas in the form
\begin{equation} \label{pre_one}
p(\vr,\vt) = (\gamma-1) \vr e(\vr,\vt).
\end{equation}
The monoatomic gas corresponds to the value $\gamma=\frac 53$.
Recall that the monoatomic gas law was carefully studied in book \cite{FNbook}, including detailed computations in this case, the generalization for arbitrary $\gamma >1$ can be found in \cite{NP1}. Below, we state the basic points.

We postulate the validity of the Gibbs relation
\begin{equation} \label{Gibbs}
\frac{1}{\vartheta} \Big(De(\vr,\vt) + p(\vr,\vt)D\Big(\frac{1}{\vr}\Big)\Big) = Ds(\vr,\vt)
\end{equation}
which is equivalent to the Maxwell relation
\begin{equation} \label{Maxwell}
\pder{e(\vr,\vt)}{\vr} = \frac{1}{\vr^2} \Big(p(\vr,\vt) - \vt
\pder{p(\vr,\vt)}{\vt}\Big).
\end{equation}
Above, Gibbs relation \eqref{Gibbs} defines a new quantity, the specific entropy, up to an additive constant. It formally satisfies (on the level of smooth solutions)
\begin{equation}\label{entropy_equality}
\Div(\vr s\vu) +\Div\Big(\frac{\vc{ q}}\vt\Big)=\frac
{\tn{S}:\nabla\vu}\vt-\frac {\vc{ q}\cdot\nabla\vt}{\vt^2}.
\end{equation}

Furthermore, due to \eqref{pre_one} and \eqref{Maxwell}, if $p \in C^1((0,\infty)^2)$, then it has necessarily
the form
\begin{equation} \label{pre_two}
p(\vr,\vt) = \vt^{\frac{\gamma}{\gamma-1}} P\Big(\frac{\rho}{\vt^{\frac {1}{\gamma-1}}}\Big),
\end{equation}
where $P \in C^1((0,\infty))$.

We will assume that
\begin{equation} \label{pre_three}
\begin{array}{c}
\displaystyle P(\cdot) \in C^1([0,\infty)) \cap C^2((0,\infty)), \\
\displaystyle P(0) = 0, \quad P'(0) = p_0 >0, \quad P'(Z) >0, \quad Z>0,\\
\displaystyle \lim_{Z \to \infty} \frac{P(Z)}{Z^{\gamma}} = p_\infty >0, \\[10pt]
\displaystyle 0 < \frac{1}{\gamma-1} \frac{\gamma P(Z) - Z P'(Z)}{Z} \leq C_5 <\infty, \quad Z>0.
\end{array}
\end{equation}
The form of the pressure together with these assumptions leads to the following properties of the pressure $p$ and the specific internal energy $e$.

For $K_0$ a fixed constant
\begin{equation} \label{pre_prop1}
\begin{array}{lclcl}
\displaystyle
C_6 \vr \vt &\leq & p(\vr,\vt) & \leq & C_7 \vr\vt, \mbox{ for } \vr \leq K_0 \vt^{\frac{1}{\gamma-1}}, \\
C_8 \vr^{\gamma} &\leq & p(\vr,\vt) & \leq & C_9 \left\{\begin{array}{ll}
\vt^{\frac{\gamma}{\gamma-1}}, &  \mbox{ for } \vr \leq K_0 \vt^{\frac{1}{\gamma-1}}, \\
\vr^\gamma, & \mbox{ for } \vr > K_0 \vt^{\frac{1}{\gamma-1}}.
\end{array}\right.
\end{array}
\end{equation}
Further
\begin{equation} \label{pre_prop2}
\begin{array}{c}
\displaystyle \pder{p(\vr,\vt)}{\vr} > 0 \qquad \mbox{ in } (0,\infty)^2, \\
\displaystyle p = d \vr^\gamma +p_m(\vr,\vt) , \quad d>0, \qquad \mbox{ with } \quad \pder
{p_m(\vr,\vt)}{\vr} \geq 0 \qquad \mbox { in } (0,\infty)^2.
\end{array}
\end{equation}
The specific internal energy defined by \eqref{pre_one} fulfills
\begin{equation} \label{ene_prop}
\left.\begin{array}{c}
\displaystyle \frac{1}{\gamma-1} p_\infty \vr^{\gamma-1} \leq e(\vr,\vt) \leq C_{10} (\vr^{\gamma-1} + \vt),  \\
\displaystyle
\pder{e(\vr,\vt)}{\vr} \vr \leq C_{10} (\vr^{\gamma-1} + \vt)
\end{array}\right\} \mbox{ in } (0,\infty)^2.
\end{equation}
Then the specific entropy has the form
\begin{equation} \label{ent_one}
s(\vr,\vt) = S\Big(\frac{\vr}{\vt^{\frac{1}{\gamma-1}}}\Big), \qquad S'(Z) = -\frac{1}{\gamma-1} \frac{\gamma P(Z)-ZP'(Z)}{Z^2}<0,
\end{equation}
as well as
\begin{equation} \label{ent_two}
\begin{array}{c}
\displaystyle \pder{s(\vr,\vt)}{\vr} = \frac{1}{\vt}\Big(-\frac{p(\vr,\vt)}{\vr^2} +
\pder{e(\vr,\vt)}{\vr}\Big) = -\frac{1}{\vr^2} \pder{p(\vr,\vt)}{\vt}, \\
\displaystyle \pder{s(\vr,\vt)}{\vartheta} = \frac{1}{\vartheta} \pder{e(\vr,\vt)}{\vartheta}
= \frac{1}{\gamma-1} \frac{\vt^{\frac{1}{\gamma-1}}}{\vr}
\Big(\gamma P\Big(\frac{\vr}{\vt^{\frac{1}{\gamma-1}}}\Big) - \frac{\vr}{\vt^{\frac{1}{\gamma-1}}}
P'\Big(\frac{\vr}{\vt^{\frac{1}{\gamma-1}}}\Big)\Big) >0.
\end{array}
\end{equation}
We also have for suitable choice of the additive constant in the definition of the specific entropy
\begin{equation} \label{ent_three}
\begin{array}{lclcl}
\displaystyle |s(\vr,\vt)| & \leq & C_{11}(1+ |\ln \vr| + |\ln \vt|) \qquad & \mbox{ in } & (0,\infty)^2, \\
\displaystyle
|s(\vr,\vt)| & \leq & C_{12} (1 + |\ln \vr|)  \qquad & \mbox{ in } & (0,\infty) \times (1,\infty), \\
\displaystyle s(\vr,\vt) &\geq & C_{13} >0  \qquad & \mbox{ in } & (0,1) \times (1,\infty), \\
\displaystyle s(\vr,\vt) & \geq & C_{14}(1+ \ln \vt)  \qquad & \mbox{ in } & (0,1) \times (0,1).
\end{array}
\end{equation}
Note that due to \eqref{ent_one} this choice of the thermodynamic potentials can be compatible with the Third law of thermodynamics which states that the entropy should go to zero when the temperature goes to zero. Indeed, the specific entropy has  a limit for $\vt\to 0^+$ which can be either equal to $-\infty$ or to a finite number. If we postulate that the limit is finite, by a suitable choice of the additive constant we may achieve that 
\begin{equation} \label{ent_ThirdLaw}
\lim_{\vt\to 0^+} s(\vr,\vt) = 0.
\end{equation}
Then as a consequence we may replace \eqref{ent_three} by
\begin{equation} \label{ent_four}
\begin{aligned}
|s(\vr,\vt)| & \leq  C_{11}(1+ |\ln \vr| + [\ln \vt]^+) \qquad \mbox{ in }  (0,\infty)^2, \\
s(\vr,\vt) &\geq 0 \qquad \mbox{ in }  (0,\infty)^2.
\end{aligned}
\end{equation}
This observation will be important in the part devoted to the Dirichlet boundary conditions for the temperature. On the other hand, in case of the heat flux boundary condition, it is possible to replace the complicated form of the pressure presented above by a simplified version with similar asymptotic behaviour
\begin{equation} \label{pre_four}
p(\vr,\vt) = d \vr^\gamma + a \vr\vt,
\end{equation}
where $a$ and $d$ are positive constants. This due to \eqref{Maxwell} implies the form of the internal energy 
\begin{equation} \label{ene_two}
e(\vr,\vt) = e_0(\vt) + d \frac{\vr^{\gamma-1}}{\gamma-1},
\end{equation}
where the function of the temperature $e_0$ is often assumed to be linear, $e_0(\vt) = b\vt$. This leads to the form of the entropy
\begin{equation} \label{ent_five}
s(\vr,\vt) = b \ln \vt - a \ln \vr + a_0
\end{equation}
which is clearly not compatible with the Third law of thermodynamics.

\subsection{Weak solution}
\label{2.2}

We now introduce the weak formulation of our problem \eqref{NSF_1}--\eqref{tot_mass}. We follow the strategy developed in \cite{NP1} or \cite{P}.  We assume that all the functions which appear in the integrals are such that all integrals are finite and we specify the function spaces later. 

The weak formulation of the continuity equation is independent of the velocity boundary conditions,
\begin{equation} \label{cont_weak}
\int_\Omega \vr\vu \cdot \nabla \psi \dx = 0 \qquad \forall \psi \in C^1(\overline{\Omega}).
\end{equation}
An equivalent possibility is to assume the integration over $\R^3$ with $\psi \in C^1_0(\R^3)$ (compactly supported differentiable functions in $\R^3$), where the density is extended by zero and the velocity in such a way that it remains in the $W^{1,p}(\Omega)$ space; in particular, for homogeneous Dirichlet boundary conditions for the velocity, the extension is by zero. In the following analysis, the renormalized formulation of the continuity equation will play an important role. This means
\begin{equation} \label{cont_renor_weak}
\int_\Omega \Big(b(\vr) \vu \cdot \nabla \psi + (b(\vr)-b'(\vr)\vr)\Div \vu \psi\Big) \dx = 0 \qquad \forall \psi \in C^1(\overline{\Omega})
\end{equation}
for any $b\in C^1([0,\infty))$ such that $b(0)=0$, $b'(z) = 0$ for $z\geq M$ with some $M>0$. Based on the integrability of both $\vr$ and $\vu$, the validity can be extended for a bigger class of functions $b$. Again, using the same extension as above outside $\Omega$, the integration over $\Omega$ can be replaced by the integration over $\R^3$ with test functions from $C^1_0(\R^3)$.

Next we consider the momentum equation. Here the situation differs for different boundary conditions for the velocity. In the case of homogeneous Dirichlet boundary conditions we consider
\begin{equation} \label{mom_weak_Dir}
\int_\Omega \Big[\big(\vr \vu \otimes \vu -\tn{S}\big):\nabla \vcg{\varphi} + p \Div \vcg{\varphi}\Big]\dx + \int_\Omega \vr\vc{f}\cdot \vcg{\varphi}\dx=0
\end{equation}
for any $\vcg{\varphi} \in C^1_0(\Omega)$.  The situation for the Navier (slip) boundary conditions is slightly different. Here the corresponding formulation is
\begin{equation} \label{mom_weak_Nav}
\int_\Omega \Big[\big(\vr \vu \otimes \vu -\tn{S}\big):\nabla \vcg{\varphi} + p \Div \vcg{\varphi}\Big]\dx -\sigma \int_{\partial \Omega} \vu \cdot \vcg{\varphi} \dS  + \int_\Omega \vr\vc{f}\cdot \vcg{\varphi}\dx=0
\end{equation}  
for any $\vcg{\varphi} \in C^1(\overline{\Omega})$ such that $\vcg{\varphi} \cdot \vc{n} =0$ at $\partial \Omega$.

Next we consider the balance of energy. Here the situation is more complex. First, we present the situation with the heat flux boundary conditions \eqref{temp_heatflux}. 

The weak formulation of the total energy balance reads as
\begin{equation} \label{ene_weak}
-\int_\Omega \Big(\frac 12 \vr |\vu|^2 \vu + \vr e \vu +\vc{q} -\tn{S} \vu + p\vu\Big)\cdot \nabla \psi \dx + \int_{\partial \Omega} \Big(L(\vt-\vt_0) + \sigma |\vu|^2 \Big)\psi \dS = \int_\Omega \vr \vc{f} \cdot \vu \psi \dx 
\end{equation}
for all $\psi \in C^1(\overline{\Omega})$. Above, if we formally set $\sigma=0$, we also include the case of the Dirichlet boundary conditions for the velocity.
However, in some situations the terms $\vr|\vu|^2 \vu$, $\tn{S}\vu$ or $p\vu$ might be not integrable. In these situations we replace the total energy balance by the total energy equality (i.e., \eqref{ene_weak} with $\psi =1$)
\begin{equation} \label{ene_equ}
\int_{\partial \Omega} \big(L(\vt-\vt_0)+ \sigma |\vu|^2\big)\dS = \int_\Omega \vr \vc{f} \cdot \vu  \dx
\end{equation}
and the entropy inequality (the equality can be obtained only in case of sufficiently smooth solutions which by nowadays available technique requires certain smallness of the data)
\begin{equation} \label{ent_ine}
\int_\Omega \Big(\frac{\tn{S}:\nabla \vu}{\vt} + \frac{\kappa(\vt) |\nabla \vt|^2}{\vt^2}\Big)\psi \dx -\int_{\partial\Omega} \frac{L(\vt-\vt_0)}{\vt}\psi\dS \leq \int_\Omega \Big(\frac{\kappa(\vt) \nabla \vt \cdot \nabla \psi}{\vt} - \vr s\vu \cdot \nabla \psi \Big)\dx  
\end{equation}
for all $\psi \in C^1(\overline{\Omega})$ such that $\psi \geq 0$ in $\overline{\Omega}$.

The situation for the Dirichlet boundary conditions for the temperature is more complex. Here we cannot replace the terms $\vc{q}\cdot \vc{n}\psi$ in the total energy balance and $\frac{\vc{q}\cdot\vc{n}}{\vt}\psi$ in the entropy inequality by the corresponding boundary conditions for the heat flux, as only the temperature is prescribed. We use the ideas from \cite{FC} adopted to the steady case in \cite{P}. We can consider the total energy balance with compactly supported test functions (if all terms are integrable), the formulation is the same for both velocity boundary conditions  
\begin{equation} \label{ene_weak_2}
\int_\Omega \Big(\frac 12 \vr |\vu|^2 \vu + \vr e \vu +\vc{q} -\tn{S} \vu + p\vu\Big)\cdot \nabla \psi \dx +  \int_\Omega \vr \vc{f} \cdot \vu \psi \dx =0 \qquad \forall \psi \in C^1_0(\Omega)
\end{equation}
and the same with the entropy inequality
\begin{equation} \label{ent_ine_2}
\int_\Omega \Big(\frac{\tn{S}:\nabla \vu}{\vt} + \frac{\kappa(\vt)|\nabla \vt|^2}{\vt^2}\Big)\psi \dx  \leq  \int_\Omega \Big(\frac{\kappa(\vt) \nabla \vt \cdot \nabla \psi}{\vt} - \vr s\vu \cdot \nabla \psi \Big)\dx \qquad \forall \psi \in C^1_0(\Omega), \, \psi \geq 0.  
\end{equation}
To have some information up to the boundary (which will be important later on, in order to get the estimates), we take the total energy balance \eqref{ene_weak} with $\psi \equiv 1$ and sum it with the entropy inequality  \eqref{ent_ine} with the test function $\psi_D$, where $\psi_D$ is the harmonic extension of the boundary data $\vt_D$. In both cases, we use the terms $\vc{q}\cdot\vc{n}$ and $\frac{\vc{q}\cdot\vc{n}}{\vt}$ instead of the terms including $L(\vt-\vt_0)$. Thus the terms $\int_{\partial \Omega} \vc{q}\cdot \vc{n}\dS$ and $\int_{\partial \Omega} \frac{\vc{q}\cdot \vc{n}}{\vt}\psi_D\dS$ cancel each other as $\psi_D=\vt$ at $\partial \Omega$ and we obtain the so called ballistic energy inequality
\begin{equation}\label{ball_ene}
\int_\Omega \Big(\frac{\tn{S}:\nabla \vu}{\vt} + \frac{\kappa(\vt)|\nabla \vt|^2}{\vt^2}\Big)\psi_D \dx + \sigma \int_{\partial \Omega} |\vu|^2\,\dS  \leq \int_\Omega \Big(\vr \vc{f} \cdot \vc{u} - \vr s\vu \cdot \nabla \psi_D + \frac{\kappa(\vt) \nabla \vt \cdot \nabla \psi_D}{\vt}\Big)\dx. 
\end{equation}
In what follows, in order to simplify the situation, we assume that the friction coefficient $\sigma=0$ and  that the domain $\Omega$ is not axially symmetric, to have the Korn inequality. The case $\sigma >0$ and the domain not axially symmetric leads to exactly the same result, if the domain is axially symmetric, then it leads to certain technicalities which we prefer to avoid (cf. \cite{JNP}).

We can now turn our attention to the weak formulations of our problems.	We start with the heat flux boundary conditions. Recall that we set $\sigma =0$, so the weak formulations for the momentum differ only by the fact whether $\vu$ and the test functions $\vcg{\varphi}$ are zero at the boundary or whether only their normal parts are zero (in the corresponding sense). In what follows, we assume that all functions in the corresponding definition are such that all integrals are finite and we assume validity of the constitutive relations from Section \ref{2.1}.

\begin{df}[Heat flux boundary condition for the temperature] \label{d1}
We say that the triple $(\vr,\vu,\vt)$ is a renormalized variational entropy solution to problem \eqref{NSF_1}, \eqref{vel_Dir} (or \eqref{vel_Nav} with $\sigma=0$), \eqref{temp_heatflux} and \eqref{tot_mass} provided $\vr \in L^{\gamma +\Theta}(\Omega)$ for some $\Theta >0$, $\vu \in W^{1,p}_0(\Omega)$ (or $\vu \in W^{1,p}(\Omega)$, $\vu \cdot \vc{n}=0$ at $\partial \Omega$), $\vt \in W^{1,r}(\Omega) \cap L^{3m}(\Omega)$ for some $r>1$, $p=\frac{6m}{3m-1+\alpha}$ and we have \eqref{cont_weak}, \eqref{cont_renor_weak}, \eqref{mom_weak_Dir} (or \eqref{mom_weak_Nav} with $\sigma =0$), \eqref{ene_equ} as well as \eqref{ent_ine}.

If, moreover, we also have \eqref{ene_weak}, the solution is called a renormalized weak solution.  
\end{df} 

Next we consider the Dirichlet boundary condition for the temperature.

\begin{df}[Dirichlet condition for the temperature] \label{d2}
We say that the triple $(\vr,\vu,\vt)$ is a renormalized ballistic energy variational entropy solution to problem \eqref{NSF_1}, \eqref{vel_Dir} (or \eqref{vel_Nav} with $\sigma=0$), \eqref{temp_Dir} and \eqref{tot_mass} provided $\vr \in L^{\gamma +\Theta}(\Omega)$ for some $\Theta >0$, $\vu \in W^{1,p}_0(\Omega)$ (or $\vu \in W^{1,p}(\Omega)$, $\vu \cdot \vc{n}=0$ at $\partial \Omega$), $\vt \in W^{1,r}(\Omega) \cap L^{3m}(\Omega)$ for some $r>1$, $p=\frac{6m}{3m-1+\alpha}$, $\vt=\vt_D$ in the sense of traces at $\partial \Omega$ and we have \eqref{cont_weak}, \eqref{cont_renor_weak}, \eqref{mom_weak_Dir} (or \eqref{mom_weak_Nav} with $\sigma =0$) \eqref{ent_ine_2} as well as \eqref{ball_ene} with $\psi_D$ the harmonic extension of $\vt_D$.

If, moreover, we also have \eqref{ene_weak_2}, the solution is called a renormalized ballistic energy weak solution.  
\end{df}

\subsection{Main results}
\label{2.3}

The main results proved in this papers are summarized below.

\begin{thm} \label{t_heat_flux}
Let $\Omega \subset \R^3$ be a $C^2$ bounded domain. Let $M>0$ in \eqref{tot_mass}, $L>0$, $\vt_0 \in L^1(\partial \Omega)$, $\vt_0 \geq T_0 >0$ a.e. at $\partial \Omega$ in \eqref{temp_heatflux} and $\vc{f}\in L^\infty(\Omega)$. Let 
\begin{equation} \label{t1_cond1}
\begin{aligned}
&\alpha \in [0,1], \qquad m >\max \Big\{\frac{1+\alpha}{3}, \frac{1-\alpha}{2}\Big\}, \\ 
& \gamma >\max \Big\{\frac{3m}{2m-1+\alpha}, \frac{3m(1+\alpha) + 2(m+1-\alpha)}{3m(1+\alpha)}, \frac{3m+1+\alpha}{3m-1+\alpha}\Big\}.
\end{aligned}
\end{equation}
Then there exists a renormalized variational entropy solution to problem \eqref{NSF_1}, \eqref{vel_Dir} (or \eqref{vel_Nav} with $\sigma=0$), \eqref{temp_heatflux} and \eqref{tot_mass} in the sense of Definition \ref{d1}. On the other hand, if
\begin{equation} \label{t1_cond2}
\alpha \in [0,1], \qquad m >1, \qquad \gamma > \frac{5m-1+\alpha}{3(m-1+\alpha)},
\end{equation}
then the constructed solution is a renormalized weak solution to the same problem. 
\end{thm}

\begin{thm} \label{t_Dirichlet_bc}
Let $\Omega \subset \R^3$ be a $C^2$ bounded domain. Let $M>0$ in \eqref{tot_mass}, $\vc{f} \in L^\infty(\Omega)$ and let the harmonic extension of $\vt_D$ to $\Omega$ be a strictly positive $W^{2,q}(\Omega)$ function for some $q>3$  in \eqref{temp_Dir}. Let either 
\begin{equation}\label{t2_cond0}
    \alpha \in [0,1], \qquad m > 1, \qquad \gamma > \frac{5m-1+\alpha}{3(m-1+\alpha)},
\end{equation}
or 
\begin{equation} \label{t2_cond1}
\begin{aligned}
&\alpha \in \Big(0,\frac 12\Big], \qquad m \in (1-\alpha,1], \\ 
& 6m\gamma^2 (m-1+\alpha) -\gamma [2(m+1-\alpha)^2 + (m-1+\alpha)(9m+1+\alpha)] \\
& \ \ \ + (m-1+\alpha)(3m+1+\alpha) >0,
\end{aligned}
\end{equation}
or
\begin{equation} \label{t2_cond2}
\begin{aligned}
&\alpha \in \Big(\frac 12,1\Big], \qquad m \in \Big(\frac{1+\alpha}{3},1\Big], \\ 
& 6m\gamma^2 (m-1+\alpha) -\gamma [2(m+1-\alpha)^2 + (m-1+\alpha)(9m+1+\alpha)] \\
& \ \ \ \ + (m-1+\alpha)(3m+1+\alpha) >0, \\
&\gamma > \frac{3m+1+\alpha}{3m-1+\alpha}.
\end{aligned}
\end{equation}
Then there exists a ballistic energy variational entropy solution to problem \eqref{NSF_1}, \eqref{vel_Dir} (or \eqref{vel_Nav} with $\sigma=0$), \eqref{temp_Dir} and \eqref{tot_mass} in the sense of Definition \ref{d2}. Moreover, if \eqref{t2_cond0} is satisfied then the constructed solution is also a ballistic energy weak solution.
\end{thm}

\begin{rmk} \label{r_t2}
Assuming $\vt_D \in W^{2-\frac 2q,q}(\partial\Omega)$ and for the continuous representative $0<C_1 \leq \vt_D\leq C_2<\infty$ at $\partial \Omega$, the assumptions on the harmonic extension of $\vt_D$ from the theorem above are clearly fulfilled.
\end{rmk}

The results will be proved in Sections \ref{4} and \ref{5}. Before doing so, let us look closer at the situations above for a particular choice of $\alpha$.

{\bf Case $\alpha=1$}:
Indeed, in this case we simply reobtain the results from \cite{NP1} (the case of Theorem \ref{t_heat_flux}) or from \cite{P} (the case of Theorem \ref{t_Dirichlet_bc}). Note that it does not follow immediately, it is necessary to check that several conditions are less restrictive and can be skipped. More precisely, it leads for the heat flux boundary conditions to $m>\frac 23$ and $\gamma > \max\{\frac 32,\frac{3m+2}{2}\}$ for the variational entropy solutions and to $m>1$ and $\gamma >\frac 53$ for the weak solutions. The case of Dirichlet boundary conditions for the temperature is more complex, however also here the conditions reduce to conditions described in \cite{P}.

\medskip

{\bf Case $\alpha=\frac 12$}:
This dependence is sometimes used to model the dependence of the viscosity on the temperature. In the case of heat flux boundary conditions and variational entropy solution, the conditions reduce to
$$
m> \frac 12, \qquad \gamma > \frac{6m}{4m-1}
$$
while for the weak solution we need 
$$
m>1, \qquad \gamma > \frac{10m-1}{3(2m-1)}.
$$

For the Dirichlet boundary condition for the temperature the set of restrictions is more complex. We require for the ballistic energy variational solution
$$
\begin{aligned}
&m>\frac 12, \qquad \gamma > \frac{10m-1}{3(2m-1)} \\
&12m \gamma^2 (2m-1) -\gamma(44m^2-4m-1) + 3(4m^2-1)>0.
\end{aligned}
$$
For ballistic energy weak solution we have conditions
$$
m>1, \qquad \gamma > \frac{10m-1}{3(2m-1)}
$$

\medskip

{\bf Case $\alpha=0$}:
The situation with heat flux boundary conditions for the temperature is interesting in comparison with the results from \cite{MP1} or \cite{MP2}. Note that in the latter which deals with both the Dirichlet  and Navier boundary conditions for the velocity, the restrictions were $\gamma >\frac 73$ and $m> \frac{3\gamma-1}{3\gamma-7}$, or, equivalently, $m>1$ and $\gamma > \frac{7m-1}{3(m-1)}$. The results obtained in this paper allow much larger set of $m's$ and $\gamma's$, namely for the variational entropy solution we require only
$$
m >\frac 12, \qquad \gamma > \max\Big\{\frac{3m}{2m-1}, \frac{5m+2}{3m}\Big\}
$$
(note that for $m\in (\frac 12,2)$ the relevant restrictions is $\gamma > \frac{3m}{2m-1}$, while for $m>2$ we require $\gamma > \frac{5m+2}{3m}$) and for the weak solution we need
$$
m >1, \qquad \gamma >\frac{5m-1}{3(m-1)}.
$$

Finally, for the Dirichlet boundary conditions for the temperature we require for the ballistic energy variational entropy solutions
$$ 
m>1, \qquad \gamma > \frac{5m-1}{3(m-1)}
$$
and the solution is in fact in this range also ballistic energy weak solution.

\section{Preliminaries}
\label{3}

In this section we recall the basic results needed in the next two sections. First note that
$$
\Big(\nabla \vu + \nabla \vu^T -\frac 23 \Div \vu \tn{I}\Big): \nabla \vu \geq \frac 12 \Big(\nabla \vu + \nabla \vu^T -\frac 23 \Div \vu\tn{I}\Big): \Big(\nabla \vu + \nabla \vu^T -\frac 23 \Div \vu\tn{I}\Big);
$$
whence
$$
\int_\Omega \Big|\Big(\nabla \vu + \nabla \vu^T -\frac 23 \Div \vu \tn{I}\Big): \nabla \vu\Big|^{\frac p2}\dx \geq C \int_\Omega \Big|\nabla \vu + \nabla \vu^T -\frac 23 \Div \vu \tn{I}\Big|^p \dx.
$$
By employing \cite[Theorem 11.22]{FNbook} on the extension by zero of $\vu$ we have
\begin{lem}[Korn inequality I] \label{Korn_1}
Let $\Omega \in C^{0,1}$. Let $1<p<\infty$ and let the tensor function $\tn{S}$ satisfy \eqref{stress}. Then there exists a constant $C>0$ such that
$$
\int_\Omega \Big|\Big(\nabla \vu + \nabla \vu^T -\frac 23 \Div \vu \tn{I}\Big): \nabla \vu\Big|^{\frac p2}\dx \geq C\|\vu\|_{1,p}^p
$$
for any $\vu \in W^{1,p}_0(\Omega)$.
\end{lem}

Now, repeating the steps of the proof of \cite[Lemma 2.3]{JNP} (the proof works for any $1<p<\infty$ exactly as for $p=2$ for which the proof is performed therein) we get
\begin{lem}[Korn inequality II] \label{Korn_2}
Let $\Omega \in C^{0,1}$ be not axisymmetric. Let $1<p<\infty$ and let the tensor function $\tn{S}$ satisfy \eqref{stress}. Then there exists a constant $C>0$ such that
$$
\int_\Omega \Big|\Big(\nabla \vu + \nabla \vu^T -\frac 23 \Div \vu \tn{I}\Big): \nabla \vu\Big|^{\frac p2}\dx \geq C\|\vu\|_{1,p}^p
$$
for any $\vu \in W^{1,p}_{\vc n}(\Omega):= \{\vu \in W^{1,p}(\Omega): \vu\cdot \vc{n} = 0 \text{ in the sense of traces at }\partial\Omega\}$. 
\end{lem} 

Next we deal with a few technical results which play an important role in the proof of strong convergence for the density sequence. For $v$ a scalar function we denote
\begin{equation} \label{Riesz_1}
({\cal R}[v])_{ij} = ((\nabla \otimes \nabla)\Delta^{-1})_{ij} v = {\cal F}^{-1}\Big[\frac{\xi_i \xi_j}{|\xi|^2}{\cal F}(v)(\xi)\Big],
\end{equation}
and for $\vc{u}$ a vector valued function
\begin{equation} \label{Riesz_2}
({\cal R}[\vc{u}])_{i} = \sum_{j=1}^3((\nabla \otimes \nabla)\Delta^{-1})_{ij} u_j = \sum_{j=1}^3 {\cal F}^{-1}\Big[\frac{\xi_i \xi_j}{|\xi|^2}{\cal F}(u_j)(\xi)\Big],
\end{equation} 
with ${\cal F}(\cdot)$ the Fourier transform. We have (see \cite[Theorems 11.32--11.35]{FNbook})

\begin{lem}[Commutators I] \label{Comm_1}
\medskip
Let $\vc{u}_n \rightharpoonup \vc{u}$ in
$L^p(\R^3)$, $v_n \rightharpoonup v$ in
$L^q(\R^3)$, where
$$
\frac 1p + \frac 1q = \frac 1s <1.
$$
Then
$$
v_n {\cal R}[\vc{u}_n] - {\cal
R}[{v}_n] \vc{u}_n \rightharpoonup v {\cal
R}[\vc{u}] - {\cal R}[{v}] \vc{u}
$$
in $L^s(\R^3)$.
\end{lem}

\begin{lem}[Commutators II] \label{Comm_2}
\medskip
Let $w \in W^{1,r}(\R^3)$, $\vc{z} \in L^p(\R^3)$, $1<r<3$,
$1<p<\infty$, $\frac 1r + \frac 1p -\frac 13 <\frac 1s <1$. Then
for all such $s$ we have
$$
\|{\cal R}[w\vc{z}] - w{\cal R}[\vc{z}] \|_{a,s,\Rr^3} \leq C
\|w\|_{1,r,\Rr^3} \|\vc{z}\|_{p,\Rr^3},
$$
where $\frac a 3 = \frac 1s + \frac 13 -\frac 1p - \frac 1r$. Here, $\|\cdot\|_{a,s, \Rr^3}$ denotes the norm in
the Sobolev--Slobodetskii space $W^{a,s}(\R^3)$.
\end{lem}

The following result is classical and its proof can be found, e.g., in \cite[Lemma 11.12]{FNbook}.

\begin{lem}[Friedrich commutator lemma] \label{Fried}
Let $N\geq 2$, $q$, $r \in (1,\infty)$, $\frac{1}{r} + \frac{1}{\beta} =  \frac 1q \in (0,1]$. Let
$$
\vr \in L_{\rm loc}^\beta(\Omega), \qquad \vu \in W^{1,r}_{\rm loc}(\Omega).
$$
Then
$$
\Div(S_\varepsilon[\vr\vu]) - \Div (S_\varepsilon[\vr]\vu) \to 0
$$
for $\varepsilon \to 0^+$ in $L^q_{\rm loc}(\Omega)$, where $S_\varepsilon$ is the molifier in spatial variables. 
\end{lem}

\begin{lem}[Weak convergence and monotone operators] \label{Monote_Weak_limit}
Let the couple of nondecreasing functions $(P,G)$ be in $C(\R) \times C(\R)$. 
  Assume that $\{\vr_n\}_{n=1}^\infty\subset
L^1(\Omega)$ is a sequence such that
$$
\left.\begin{array}{c}
P(\vr_n) \rightharpoonup \overline{P(\vr)} \\
G(\vr_n) \rightharpoonup \overline{G(\vr)} \\
P(\vr_n)G(\vr_n) \rightharpoonup \overline{P(\vr)G(\vr)}
\end{array} \right\} \mbox{ in } L^1(\Omega).
$$
Then
$$
\overline{P(\vr)}\, \, \overline{G(\vr)} \leq
\overline{P(\vr)G(\vr)}
$$
a.e. in $\Omega$, where the bar over a function, e.g. $\overline P(\vr)$ denotes the weak limit of the sequence $P(\vr_n)$.
\end{lem}

The proof can be found e.g. in \cite[Theorem 11.26]{FNbook}

\begin{lem}[Weak convergence and convex functions] \label{weak_convex}
Let $\{\vr_n\}_{n=1}^\infty$ be a sequence of functions in $L^1(\Omega)$ such that
$$
\vr_n \rightharpoonup \vr \qquad \text{ in } L^1(\Omega).
$$
Let $\Phi$ be a continuous convex function. Then
$$
\int_\Omega \Phi(\vr)\dx \leq \liminf_{n\to \infty} \int_\Omega \Phi(\vrn) \dx.
$$
Moreover, if $\Phi(\vrn) \rightharpoonup \overline{\Phi(\vr)}$ weakly in  $L^1(\Omega)$, then
$$
\Phi(\vr) \leq \overline{\Phi(\vr)} \qquad \text{ a.a. in } \Omega.
$$
If, in addition, $\Phi$ is strictly convex on an open convex set $U\subset \R^N$ and
$$
\Phi(\vr) = \overline{\Phi(\vr)} \qquad \text{ a.a. in } \Omega,
$$
then
$$
\vrn(x) \to \vr(x) \quad \text{for a.e. } x\in \{y\in \Omega: \vr(x) \subset U\}.
$$
\end{lem}

The proof can be found e.g. in \cite[Theorem 11.27]{FNbook}

For the estimates of the density we need a special branch of solutions to the following problem
\begin{equation} \label{Bog}
\begin{array}{c}
\Div \vcg{\Phi} = f \quad \mbox{ in } \Omega,\\
\vcg{\Phi} = \vc{0}\quad \mbox{ at } \partial \Omega.
\end{array}
\end{equation}

We have (see, e.g., \cite[Theorem 11.17]{FNbook})
\begin{lem} \label{Bogovskii}
Let $f \in L^p(\Omega)$, $1<p<\infty$, $\intO{f} = 0$, $\Omega \in C^{0,1}$. Then there exists a solution to \eqref{Bog} and a constant $C>0$ independent of $f$ such that
\begin{equation} \label{2.6}
\|\vcg\Phi\|_{1,p} \leq C \|f\|_p.
\end{equation}
Moreover, the solution operator ${\cal B}$: $f \mapsto \vcg{\Phi}$ is linear.
\end{lem}


\section{Proof for the heat flux boundary condition}
\label{4}

We will skip the construction of the approximate solutions since we can employ, with minor changes the construction from \cite{NP1}. The only changes are connected with the different form of the viscosities and this does not lead to any big differences in the procedure. We will therefore present only the weak compactness of the solutions; more precisely, assuming a sequence of solutions in the sense of Definition \ref{d1} (constructed, e.g., by a convergent sequence of the total masses $M_n$) we first show certain estimates for the triple $(\vr_n,\vu_n,\vt_n)$ and then we will verify that they guarantee that the corresponding limit functions fulfill the same problem with the limit total mass $M=\lim_{n\to \infty} M_n$. In fact, we will be more detailed in the part of estimates while for the study of the convergence we will only comment carefully the parts which are different with respect to \cite{NP1} and only sketch the parts which are exactly the same. A similar approach will be employed  also for the next Section \ref{5} dealing with the Dirichlet boundary conditions for the temperature. 

\subsection{Energy and entropy estimates} \label{4.1}

Let us consider a sequence of solutions $(\vrn,\vun,\vtn)$ to our problem \eqref{NSF_1}--\eqref{temp_heatflux}, \eqref{tot_mass} (either with Dirichlet or Navier boundary conditions for the velocity) with constitutive relations satisfying the assumptions from Subsection \ref{2.1}. In fact, we can work in this section either with the pressure form \eqref{pre_two} and the related entropy \eqref{ent_one} with or without the Third law of thermodynamics \eqref{ent_ThirdLaw}, or with \eqref{pre_four}, \eqref{ene_two} and \eqref{ent_five}.

Assume that the sequence satisfies Definition \ref{d1} (with the corresponding form for the momentum equation related with the given velocity boundary conditions) and assume that the solution is a renormalized weak one and sufficiently integrable so that in particular the test function based on Bogovskii operator in the momentum equation is allowed.

The entropy inequality with the test function $\psi =1$ yields
\begin{equation} \label{ent_1}
\int_\Omega \Big(\frac{\tn{S}(\vtn,\nabla \vun):\nabla \vun}{\vtn} + \frac{\kappa(\vtn)|\nabla \vtn|^2}{\vtn^2}\Big)\,\dx + L \int_{\partial \Omega}\frac{\vt_0}{\vtn} \,\dS = L|\partial \Omega| 
\end{equation}
and the total energy balance with the test function $\psi=1$ yields
\begin{equation} \label{ene_1}
L\int_{\partial \Omega} \vtn \,\dS = \int_\Omega \vrn \vu_n \cdot \vc{f} \,\dx + L \int_{\partial \Omega} \vt_0 \,\dS.
\end{equation} 
From \eqref{ent_1} we get by H\"older's inequality for $p= \frac{6m}{3m+1-\alpha}$ and by the corresponding form of the Korn's inequality (Lemma \ref{Korn_1} or \ref{Korn_2})
\begin{equation} \label{48a}
\begin{aligned}
\|\vun\|_{1,p}^p &\leq C   \Big(\int_\Omega \frac{\tn{S}(\vtn,\nabla\vun):\nabla \vun}{\vtn}\, \dx\Big)^{\frac{p}{2}} \Big(\int_\Omega \vtn^{\frac{1-\alpha}{2} \frac{2p}{2-p}}\,\dx \Big)^{\frac{2-p}{2}} \\
&\leq C \|\vtn\|_{3m}^{3m\frac{1-\alpha}{3m+1-\alpha}}, 
\end{aligned}  
\end{equation}
whence
\begin{equation} \label{vel_1}
\|\vun\|_{1,p} \leq C \|\vtn\|_{3m}^{\frac{1-\alpha}{2}}.
\end{equation}
Using the total energy equality \eqref{ene_1} together with estimate \eqref{vel_1} we conclude
$$
\|\vtn\|_{3m} \leq C \Big(1+ \int_\Omega |\vc{f}||\vun| \vrn \,\dx\Big) \leq C \big(1+ \|\vun\|_{\frac{3p}{3-p}}\|\vrn\|_{\frac{3p}{4p-3}}\big),
$$
since $1<p\leq 2$. Plugging this estimate into the estimate of velocity we get 
$$
\|\vun\|_{1,p} \leq C \big(1 + \|\vun\|_{\frac{3p}{3-p}} \|\vrn\|_{\frac{3p}{4p-3}}\big)^{\frac{1-\alpha}{2}} \leq \frac 12 \|\vun\|_{1,p} + C\big(1+ \|\vrn\|_{\frac{3p}{4p-3}}^{\frac{1-\alpha}{1+\alpha}}\big).
$$
Therefore we conclude
\begin{equation} \label{est_1}
\begin{aligned}
\|\vun\|_{1,p} &\leq C \big(1+ \|\vrn\|_{\frac{3p}{4p-3}}^{\frac{1-\alpha}{1+\alpha}}\big) \\
\|\vtn\|_{3m} &\leq C \big(1+ \|\vrn\|_{\frac{3p}{4p-3}}^{\frac{2}{1+\alpha}}\big)
\end{aligned}
\end{equation}
which, indeed, simplifies into the known estimates for $\alpha=1$. However, this estimate is less restrictive than a similar estimate in \cite{KNP}, where the presence of radiation lead to more restrictive estimates with respect to \eqref{est_1}.

Later on, we will interpolate the density norm between $L^1$ and $L^{\gamma(1+\Theta)}$, however, let us first explain the choice of the exponent $\gamma(1+\Theta)$.

\subsection{Estimates of the density based on Bogovskii operator} \label{4.2}

We now consider the momentum equation. For both velocity boundary conditions, the function
$$
\vcg{\varphi}:= \mathcal B \Big[p(\vrn,\vtn)^\Theta - \frac{1}{|\Omega|} \int_\Omega p(\vrn,\vtn)^\Theta \,\dx\Big]
$$
for some $\Theta>0$ is an admissible test function, recall Lemma \ref{Bogovskii}. Therefore we get
\begin{equation} \label{Bog_test}
\begin{aligned}
\int_\Omega p(\vrn,\vtn)^{1+\Theta} \,\dx &= -\int_\Omega \vrn (\vun\otimes\vun) :\nabla \vcg{\varphi} \,\dx + \int_\Omega \tn{S}(\vtn,\nabla \vun):\nabla \vcg{\varphi}\,\dx \\
&- \int_\Omega \vrn\vc{f}\cdot \vcg{\varphi}\,\dx + \int_\Omega p(\vrn,\vun)\,\dx \frac{1}{|\Omega|} \int_\Omega p(\vrn,\vtn)^\Theta\, \dx \\
& = I_1+I_2+I_3 +I_4.
\end{aligned}
\end{equation}
The left-hand side for the pressure of the form \eqref{pre_one} is bounded from below up to an multiplicative constant by
$$
\int_\Omega \big( \vrn^{\gamma(1+\Theta)} + (\vrn\vtn)^{1+\Theta} 1_{\{x\in \Omega: \vrn(x) \leq K_0 \vtn(x)^{\frac{1}{\gamma-1}}\}}\big)\,\dx.
$$
Note also that on the set $\{x\in \Omega: \vrn(x) > K_0 \vtn(x)^{\frac{1}{\gamma-1}}\}$ we have that $\vrn\vtn \leq C \vrn^\gamma$. On the other hand, for the pressure of the form \eqref{pre_four}, the lower limit is the same, only the characteristic function does not appear.

Therefore it make sense to replace in \eqref{est_1} the norm of the density in $L^{\frac{3p}{4p-3}}$ by the norm in $L^{\gamma(1+\Theta)}$. We have
\begin{equation}\label{eq:interpol}
\|\vrn\|_{\frac{3p}{4p-3}} = \|\vrn\|_{\frac{6m}{5m-1 +\alpha}} \leq \|\vrn\|_1^{1-\lambda} \|\vrn\|_{\gamma(1+\Theta)}^\lambda,
\end{equation}
where $\lambda = \frac{\gamma(1+\Theta)(m+1-\alpha)}{6m(\gamma(1+\Theta)-1)}$ provided $\frac{6m}{5m-1+\alpha} <\gamma(1+\Theta)$ (the other restriction, $\frac{6m}{5m-1+\alpha}>1$, is clearly fulfilled since $m>0$ and $\alpha \in [0,1]$).

Inserting the interpolation inequality to \eqref{est_1} and recalling that the $L^1$-norm of the density is fixed, we end up with  
\begin{equation} \label{est_2}
\begin{aligned}
\|\vun\|_{1,p} &\leq C \big(1+ \|\vrn\|_{\gamma(1+\Theta)}^{\frac{\gamma(1+\Theta)(m+1-\alpha)(1-\alpha)}{6m(\gamma(1+\Theta)-1)(1+\alpha)}}\big) \\
\|\vtn\|_{3m} &\leq C \big(1+ \|\vrn\|_{\gamma(1+\Theta)}^{\frac{\gamma(1+\Theta)(m+1-\alpha)}{3m(\gamma(1+\Theta)-1)(1+\alpha)}}\big).
\end{aligned}
\end{equation}
Recall that \eqref{est_2} holds provided $\alpha \in [0,1]$, $\gamma >1$  and $\frac{6m}{5m-1+\alpha} <\gamma(1+\Theta)$.

We now return back to the momentum equation tested by a certain Bogovskii operator and proceed in two steps. First we show that under certain restrictions on $\Theta$ the terms on the right-hand side are bounded by suitable norms of $\vrn,\vun, \vtn$ and in the second step we show, based on estimates \eqref{est_2} that we indeed obtain an estimate of the density which will result into estimates of all three functions. Before we start with it, note that the term $I_3$ is clearly of lower order with respect to, e.g., the convective term and we will not consider it. 

We first have
$$
|I_1| \leq \int_\Omega \vrn |\vun|^2 |\nabla \vcg{\varphi}| \,\dx \leq \|\nabla \vcg{\varphi}\|_{\frac{1+\Theta}{\Theta}} \|\vrn\|_{\gamma(1+\Theta)} \|\vun\|_{\frac {3p}{3-p}}^2
$$
provided
$$
\frac{2(3-p)}{3p} + \frac{1}{\gamma(1+\Theta)}+ \frac{\Theta}{1+\Theta} = \frac{m+1-\alpha}{3m} + \frac{1}{\gamma(1+\Theta)}+ \frac{\Theta}{1+\Theta} \leq 1.
$$
This yields
\begin{equation} \label{Theta_1}
\Theta  \leq \frac{\gamma(2m-1+\alpha)-3m}{\gamma(m+1-\alpha)}.
\end{equation}
Since we require $\Theta >0$, it leads to restrictions
$$
\gamma > \frac{3m}{2m-1+\alpha}, \qquad m>\frac{1-\alpha}{2}.
$$
Next
$$
\begin{aligned}
|I_2| &\leq \int_\Omega |\tn{S}(\vtn,\nabla \vun):\nabla \vcg{\varphi}|\,\dx \leq C \int_\Omega (1+\vtn)^\alpha |\nabla \vun| |\nabla \vcg{\varphi}|\,\dx \\
& \leq C (1+ \|\vtn\|_{3m}^\alpha) \|\nabla \vun\|_p \|\nabla \vcg{\varphi}\|_{\frac{1+\Theta}{\Theta}}
\end{aligned}
$$
provided
$$
\frac{\alpha}{3m}+ \frac{3m+1-\alpha}{6m}+ \frac{\Theta}{1+\Theta} \leq 1.
$$
This yields
\begin{equation} \label{Theta_2}
\Theta \leq \frac{3m-1-\alpha}{3m+1+\alpha}.
\end{equation}
The upper bound is positive, provided $m> \frac{1+\alpha}{3}$.

Finally, we consider $I_4$. We will perform the proof for the pressure of the form \eqref{pre_four}, the other situation can be estimated similarly due to lower pressure estimate presented above and \eqref{pre_prop1}. Since
$$
I_4 \leq C \int_\Omega (\vrn^\gamma + \vrn\vtn)\,\dx \int_\Omega (\vrn^\gamma + \vrn\vtn)^\Theta \,\dx,
$$
whenever we have a term with density only, by interpolation as above we may estimate this term by $\int_\Omega (\vrn^\gamma + \vrn\vtn)^{1+\Theta}\,\dx
$ with a small constant in front of it (by Young's inequality). Therefore we are left with
$$
I_{4,4}:= \int_\Omega \vrn\vtn\,\dx \int_\Omega (\vrn\vtn)^\Theta \,\dx
$$
only. We get
$$
\begin{aligned}
I_{4,4} &\leq C\Big(\int_\Omega (\vrn\vtn)^{1+\Theta}\,\dx\Big)^{\frac{1}{1+\Theta}} \Big(\int_\Omega \vtn^{3m} \,\dx\Big)^{\frac{\Theta}{3m}}
\Big(\int_\Omega \vrn^{\frac{3m}{3m-\Theta}\Theta}\Big)^{\frac{3m-\Theta}{3m}} \\
&\leq \frac 12 \int_\Omega (\vrn\vtn)^{1+\Theta}\,\dx + C \|\vtn\|_{3m}^{1+\Theta} \|\vrn^{\frac{3m\Theta}{3m-\Theta}}\|_{1}^{\frac {(3m-\Theta)(1+\Theta)}{3m\Theta}}. 
\end{aligned}
$$ 
It holds
$$
\frac{3m\Theta}{3m-\Theta} < 1.
$$
Indeed, the inequality above can be transformed into
$$
\Theta \leq  \frac{3m}{3m+1}
$$
and since $\frac{3m}{3m+1} >\frac{3m-1-\alpha}{3m+1+\alpha}$ and $\Theta \leq \frac{3m-1-\alpha}{3m+1+\alpha}$, we conclude that $\frac{3m\Theta}{3m-\Theta} < 1$.
Hence
$$
I_{4,4}\leq \frac 12 \int_\Omega (\vrn\vtn)^{1+\Theta}\,\dx + C \|\vtn\|_{3m}^{1+\Theta}.
$$
Summarizing, we proved that $I_1+I_2+I_3+I_4$ is finite for our sequence $(\vrn,\vun,\vtn)$ provided
\begin{equation} \label{cond_1}
\begin{aligned}
\gamma > \frac{3m}{2m-1+\alpha}, & \qquad m > \max \Big\{\frac{1+\alpha}{3}, \frac{1-\alpha}{2}\Big\}, \\
\Theta \leq \min \Big\{&\frac{3m-1-\alpha}{3m+1+\alpha}, \frac{\gamma(2m-1+\alpha)-3m}{\gamma(m+1-\alpha)}\Big\}.
\end{aligned}
\end{equation}
In what follows, we show that in fact under certain conditions on our parameters, we have
$$
I_1+I_2 +I_3+I_4 \leq C\Big(1+ \Big(\int_\Omega p(\vrn,\vtn)^{1+\Theta}\,\dx\Big)^\delta \Big)
$$
with some $0\leq \delta <1$ which allows to conclude the boundedness of our sequence of solutions. Here, we will use \eqref{est_2} and thus, the assumption $\frac{6m}{5m+1-\alpha} \leq \gamma(1+\Theta)$ must be also taken into account.

We have
$$
\begin{aligned}
|I_1| &\leq \|p(\vrn,\vtn)\|_{1+\Theta}^\Theta \|\vrn\|_{\gamma(1+\Theta)}\|\vun\|_{\frac{3p}{3-p}}^2 \leq \varepsilon \|p(\vrn,\vtn)\|_{1+\Theta}^{1+\Theta} + C(\varepsilon) \|\vrn\|_{\gamma(1+\Theta)}^{1+\Theta}\|\vun\|_{\frac{3p}{3-p}}^{2(1+\Theta)} \\
&\leq \varepsilon \|p(\vrn,\vtn)\|_{1+\Theta}^{1+\Theta} + C(\varepsilon) \|\vrn\|_{\gamma(1+\Theta)}^{1+\Theta} \Big(1+ \|\vrn\|_{\gamma(1+\Theta)}^{2(1+\Theta) \frac{1-\alpha}{1+\alpha} \frac{\gamma(1+\Theta)(m+1-\alpha)}{6m(\gamma(1+\Theta)-1)}}\Big).
\end{aligned}
$$
Whence, the term can be estimated, provided
$$
\gamma(1+\Theta) > 1+\Theta + 2(1+\Theta) \frac{1-\alpha}{1+\alpha} \, \frac{\gamma(1+\Theta)(m+1-\alpha)}{6m(\gamma(1+\Theta)-1)}
$$
which after standard computations leads to
$$
\Theta> \frac{(1-\alpha) \gamma (m+1-\alpha)-3m(1+\alpha) (\gamma-1)^2}{\gamma\big[3m(\gamma-1)(1+\alpha) -(1-\alpha) (m+1-\alpha)\big]}.
$$
Since $\gamma > \frac{3m}{2m-1+\alpha}$, we have $\gamma-1 > \frac{m+1-\alpha}{2m-1+\alpha}$ and then as $\frac{3m}{2m-1+\alpha}>1$,
$$
3m(\gamma-1)(1+\alpha) > \frac{3m(1+\alpha)(m+1-\alpha)}{2m-1+\alpha}> (m+1-\alpha)(1-\alpha).
$$
Thus the denominator is always positive.

Next
$$
\begin{aligned}
|I_2| &\leq (1+ \|\vtn\|_{3m}^\alpha) \|\nabla \vun\|_p \|p(\vrn,\vtn)\|_{1+\Theta}^\Theta \\
&\leq \varepsilon \|p(\vrn,\vtn)\|_{1+\Theta}^{1+\Theta} + C(\varepsilon) \Big(1+\|\vtn\|_{3m}^{\alpha(1+\Theta)}\Big) \|\nabla \vun\|_p^{1+\Theta} \\
&\leq  \varepsilon \|p(\vrn,\vtn)\|_{1+\Theta}^{1+\Theta} + C(\varepsilon) \Big(1+ \|\vrn\|_{\gamma(1+\Theta)}^{(1+\Theta)(\frac{2\alpha}{1+\alpha} +\frac{1-\alpha}{1+\alpha})\frac{\gamma(1+\Theta)(m+1-\alpha)}{6m(\gamma(1+\Theta)-1)} }\Big). 
\end{aligned}
$$
Here, the inequality
$$
\gamma(1+\Theta) >(1+\Theta) \frac{\gamma(1+\Theta)(m+1-\alpha)}{6m(\gamma(1+\Theta)-1)}
$$
is fulfilled always, since it can be rewritten as
$$
\Theta(m+1-\alpha-6m\gamma) < 6m(\gamma-1)-m-1+\alpha
$$
and the left hand side is negative while the right hand side is positive (it is enough to use $\gamma > \frac32$).

Next
$$
I_4 \leq \varepsilon \|p(\vrn,\vtn)\|_{1+\Theta}^{1+\Theta} + C \|\vtn\|_{3m}^{1+\Theta} \leq  \varepsilon \|p(\vrn,\vtn)\|_{1+\Theta}^{1+\Theta} + C(\varepsilon) \|\vrn\|_{\gamma(1+\Theta)}^{(1+\Theta) \frac{2}{1+\alpha}\frac{\gamma(1+\Theta)(m+1-\alpha)}{6m(\gamma(1+\Theta)-1)}}.
$$
We therefore require
$$
\gamma(1+\Theta) > (1+\Theta) \frac{2}{1+\alpha}\frac{\gamma(1+\Theta)(m+1-\alpha)}{6m(\gamma(1+\Theta)-1)}
$$
which results into
$$
\Theta > \frac{m+1-\alpha -3m(1+\alpha)(\gamma-1)}{3m\gamma (1+\alpha) -(m+1-\alpha)}.
$$
The denominator is always positive (as $m> \frac{1-\alpha}{2}$) while the numerator is negative. This follows due to the fact that $\gamma-1 > \frac{m+1-\alpha}{2m-1+\alpha}$ and thus 
$$
3m(1+\alpha)(\gamma-1) > \frac{3m(1+\alpha)(m+1-\alpha)}{2m-1+\alpha}>m+1-\alpha,
$$ 
since $3m(1+\alpha)>2m-1+\alpha$. Therefore the terms $I_2$ and $I_4$ do not add any restriction on $\Theta$. To conclude this part, let us verify that the condition on $\Theta$ in the interpolation of the density norm does not lead to any further restriction on our parameters. More precisely, we shall check that
\begin{equation} \label{int}
\gamma(1+ \Theta) > \frac{6m}{5m-1+\alpha}
\end{equation}
is in agreement with the upper bounds on $\Theta$.

First, let $\Theta = \frac{3m-1-\alpha}{3m+1+\alpha}$. Replacing $\Theta$ in \eqref{int} by this value leads to
\begin{equation} \label{gamma_1}
\gamma > \frac{3m+1+\alpha}{5m-1+\alpha}.
\end{equation}
However, since $\gamma > \frac{3m}{2m-1+\alpha}$, condition \eqref{gamma_1} is less restrictive provided $m>\frac 29 (\alpha+1)$ which holds true due to $m>\frac 13(\alpha+1)$.

Next, let $\Theta = \frac{\gamma(2m-1+\alpha)-3m}{\gamma(m+1-\alpha)}$. Plugging in this value into \eqref{int} we get $\frac{3m(\gamma-1)}{m+1-\alpha} > \frac{6m}{5m-1+\alpha}$ which is less restrictive since it yields $\gamma > \frac{7m+1-\alpha}{5m-1+\alpha}$
and $\frac{7m+1-\alpha}{5m-1+\alpha} < \frac{3m}{2m-1+\alpha}$ for $m > \frac{1-\alpha}{2}$.

Concluding,  the estimates
$$
\|\vrn\|_{\gamma(1+\Theta)} + \|\vtn\|_{3m} + \|\nabla (\vtn^{\frac m2})\|_2 + \|\ln\vtn\|_{1,2} + \|\vun\|_{1,p} \leq C
$$
hold true for $p = \frac{6m}{3m+1-\alpha}$, provided
$$
\begin{aligned}
\Theta&= \min\Big\{\frac{3m-1-\alpha}{3m+1+\alpha}, \frac{\gamma(2m-1+\alpha)-3m}{\gamma(m+1-\alpha)}\Big\}, \\
 \Theta &> \frac{(1-\alpha)\gamma (m+1-\alpha)-3m(1+\alpha)(\gamma-1)^2}{\gamma[3m(\gamma-1)(1+\alpha)-(1-\alpha) (m+1-\alpha)]}
\end{aligned}
$$
with extra conditions
$$
m> \max\Big\{\frac{1+\alpha}{3}, \frac{1-\alpha}{2}\Big\} \qquad \gamma > \frac{3m}{2m-1+\alpha}.
$$
In what follows, we show that also the last lower bound on $\Theta$ is not in contradiction with respect to the upper bounds on $\Theta$ in our range of parameters. 

First, let $\Theta = \frac{3m-1-\alpha}{3m+1+\alpha}$. We need to verify that under the conditions above
$$
\frac{3m-1-\alpha}{3m+1+\alpha} > \frac{(1-\alpha)\gamma (m+1-\alpha)-3m(1+\alpha)(\gamma-1)^2}{\gamma[3m(\gamma-1)(1+\alpha)-(1-\alpha) (m+1-\alpha)]}.
$$
After standard computations we end up with the quadratic inequality
$$
\gamma^2 6m (1+\alpha) -\gamma[9m(1+\alpha) + (1+\alpha)^2 +2(1-\alpha)(m+1-\alpha)] + 3m(1+\alpha) +(1+\alpha)^2 >0.
$$
It is not difficult to see that for $\gamma \geq \frac{3m}{2m-1+\alpha}$ and the lower bounds on $m$ the left-hand side of the inequality is increasing in $\gamma$. We therefore plug in the lowest value of $\gamma$ and end up with
$$
24m^3\alpha + m^2 (7+8\alpha-23\alpha^2) + m(8-22\alpha+16\alpha^2 -2\alpha^3) + (1-\alpha)^2 (1+\alpha)^2 >0
$$
which can be verified to hold. It is rather straightforward for $m\geq 1$ (it is enough to use that $m^3 \geq m^2\geq m$ and $\alpha^3 \leq \alpha^2\leq \alpha \leq 1$). For $m<1$ we need to use that $m>\frac{1+\alpha}{3}$. 

Next we take the other condition and aim to verify that
$$
\frac{\gamma(2m-1+\alpha)-3m}{\gamma(m+1-\alpha)} > \frac{(1-\alpha)\gamma (m+1-\alpha)-3m(1+\alpha)(\gamma-1)^2}{\gamma[3m(\gamma-1)(1+\alpha)-(1-\alpha) (m+1-\alpha)]}.
$$
After standard computations we end up with (we cancel $\gamma-1$ during the computations)
$$
\gamma > \frac{3m(1+\alpha)+ 2(m+1-\alpha)}{3m(1+\alpha)}
$$
which must be taken into account together with $\gamma > \frac{3m}{2m-1+\alpha}$.
Thus under the above deduced conditions we can choose sufficiently large $\Theta$ (equal to the smaller upper bound) so that the estimate on $\vrn$ holds.

Concluding, the estimates
\begin{equation} \label{est_fin}
\|\vrn\|_{\gamma(1+\Theta)} + \|\vtn\|_{3m} + \|\nabla (\vtn^{\frac m2})\|_2 + \|\ln\vtn\|_{1,2} + \|\vun\|_{1,p} \leq C
\end{equation}
hold with $p = \frac{6m}{3m+1-\alpha}$, pro\-vi\-ded
\begin{equation} \label{cond_fin}
\begin{aligned}
\Theta&= \min\Big\{\frac{3m-1-\alpha}{3m+1+\alpha}, \frac{\gamma(2m-1+\alpha)-3m}{\gamma(m+1-\alpha)}\Big\}, \\
m&> \max\Big\{\frac{1+\alpha}{3}, \frac{1-\alpha}{2}\Big\}, \qquad \gamma > \max \Big\{\frac{3m}{2m-1+\alpha}, \frac{3m(1+\alpha)+ 2(m+1-\alpha)}{3m(1+\alpha)}\Big\}.
\end{aligned}
\end{equation}
Recall that $\Theta$ is positive in this situation. Finally note that $\frac{3m}{2m-1+\alpha} > \frac{3m(1+\alpha)+ 2(m+1-\alpha)}{3m(1+\alpha)}$ for $\alpha \geq \frac 13$.

\subsection{Limit passage I}

Based on the estimates proved above we now want to pass to the limit in the weak formulations from Definition \ref{d1}. We need to do so in two steps, first we only use the estimates and thus we cannot pass to the desired limit in nonlinear functions of density which will be the goal in the next subsection. To avoid too many technicalities, we perform the proof only for the homogeneous Dirichlet boundary conditions for the velocity, the proof in the case of the Navier boundary conditions differs only in technical details connected with only normal trace of the velocity and of some test functions being zero, while for the Dirichlet ones, the full trace is zero.

Using \eqref{est_fin} we deduce existence of a subsequence (however, we relabel it so that no special notation is used, similarly also in the next section) such that
$$
\begin{aligned}
\vun \rightharpoonup \vu& \qquad \text{in }W^{1,p}_0(\Omega)\\
\vun \to \vu& \qquad \text{in }L^{r}(\Omega), \, 1\leq r < \frac {3p}{3-p} \\
\vrn \rightharpoonup \vr& \qquad \text{in }L^{\gamma(1+\Theta)}(\Omega)\\
\vtn \rightharpoonup \vt& \qquad \text{in }W^{1,r}(\Omega) \text{ for some } r>1\\
\vtn \to \vt& \qquad \text{in }L^{r}(\Omega), \, 1\leq r < 3m. \\
\end{aligned}
$$
Since also $p(\vtn,\vrn)$ is bounded in $L^{1+\Theta}(\Omega)$ for some $\Theta >0$, we know that 
$$
p(\vtn,\vrn) \rightharpoonup \overline{p(\vr,\vt)} \qquad \text{in } L^{1+\Theta}(\Omega),
$$
where we, however, do not know whether $\overline{p(\vr,\vt)} = p(\vr,\vt)$. We will use the same notation also for other nonlinear functions of density ($s(\vr,\vt)$ and $e(\vr,\vt)$). Both quantities are bounded in a better space than $L^1$ and we may pass to the limit in the weak formulations. We directly get
$$
\begin{aligned}
&\int_\Omega \vr\vu\cdot \nabla \psi \,\dx = 0 \qquad \text{for any }\  \psi \in C^1(\overline{\Omega}), \\
&\int_\Omega \Big[\big(\vr \vu\otimes\vu-\tn{S}(\vt,\nabla \vu)\big):\nabla\vcg{\varphi} + \overline{p(\vr,\vt)}\Div\vcg{\varphi}\Big]\,\dx = \int_\Omega \vr\vc{f}\cdot \vcg{\varphi}\,\dx \qquad\text{for any }\ \vcg{\varphi}\in C^1_0(\Omega) \\
&\int_\Omega \Big[\frac{\tn{S}(\vt,\vu):\nabla \vu}{\vt} + \frac{\kappa(\vt)|\nabla \vt|^2}{\vt^2}\Big]\psi \,\dx - \int_{\partial \Omega}\frac{L(\vt-\vt_0)}{\vt}\psi \,\dx \\
& \ \ \ \ \ \leq \int_\Omega\Big[\frac{\kappa(\vt)\nabla \vt }{\vt} -\vr \overline{s(\vr,\vt)} \vu\Big]\cdot \nabla \psi\,\dx
\end{aligned}
$$
for any $\psi \in C^1(\overline{\Omega})$, nonnegative.  In the last inequality, however, on the left-hand side we need to use the weak lower semicontinuity for both terms. The issue with the renormalized continuity equation is more complex and we will come to it in the next subsection. 

The limit passage in the total energy balance is more involved due to possible low integrability of some terms. We need to pass to the limit in the integrals
$$
\int_\Omega \Big(\frac 12 \vrn |\vun|^2 \vun + \vrn e(\vrn,\vtn) \vun +\vc{q}(\vrn,\vtn) -\tn{S}(\vtn,\nabla \vun) \vun + p(\vrn,\vtn)\vun\Big)\cdot \nabla \psi \dx = \sum_{i=1}^5 I_i^n.
$$
We easily get for $I_3^n$
$$
\begin{aligned}
|I_3^n| &\leq C \int_\Omega |\vc{q}(\vtn,\nabla \vtn)|\,\dx \leq C\int_\Omega |\nabla \vtn| (1+\vtn)^m\,\dx \\
&\leq C \Big(\int_\Omega |\nabla \vtn|^2 \frac{(1+\vtn)^{m}}{\vtn^2} \,\dx \Big)^{\frac 12} \Big(\int_\Omega (1+\vtn)^{m+2}\,\dx\Big)^{\frac 12}
\end{aligned}
$$
which is bounded for $m\geq 1$ since then $m+2 \leq 3m$. As the term depends linearly on $\nabla \vtn$, we can pass to the limit to get the desired integral
$$
\int_\Omega \vc{q}(\vt,\nabla \vt)\cdot \nabla \psi \,\dx
$$
provided $m>1$. 

Similarly we may proceed for $I_4^n$. We get
$$
|I_4^n| \leq C\int_\Omega |\vun||\nabla \vun| (1+\vtn)^\alpha \, \dx \leq C \|\vun\|_{\frac{3p}{3-p}} \|\nabla \vun\|_p \|1+\vtn\|_{3m}^\alpha
$$
provided ($\frac {3p}{3-p} = \frac{6m}{m+1-\alpha}$)
$$
\frac{\alpha}{3m} + \frac{3m+1-\alpha}{6m}+ \frac{m+1-\alpha}{6m} \leq 1
$$  
which leads to $m\geq 1$. Similarly, the desired limit passage can be performed provided $m>1$. 

For the kinetic energy transport term we have
$$
|I_1^n| \leq C \int_\Omega \vrn|\vun|^3\,\dx \leq C\|\vun\|_{\frac {3p}{3-p}}^3 \|\vrn\|_{\frac{p}{2p-3}}.
$$
Since $\frac{p}{2p-3} = \frac{2m}{m-1+\alpha}$, the estimate above holds provided $\gamma(1+\Theta) \geq \frac{2m}{m+1-\alpha}$. Since $1+\Theta \leq \min\{\frac{6m}{3m+1+\alpha},\frac{3m(\gamma-1)}{\gamma(m+1-\alpha)}\}$, plugging in these bounds we end up with the following restrictions on $\gamma$:
$$
\gamma \geq \frac 53, \qquad \gamma \geq \frac{3m+1+\alpha}{3(m+1-\alpha)}.
$$
Indeed, for the convergence of the integrals the inequalities must to be strict.

Last, we need to control the internal energy and the pressure terms. Since they behave similarly, we look closely only at the pressure, the internal energy follows by the same arguments. Furthermore, we concentrate on the pressure in the form \eqref{pre_four}, the other case is similar as was discussed above.

The pressure divides into two terms. First
$$
|I_{5,1}^n| \leq C \int_\Omega \vrn^\gamma |\vun|\,\dx \leq C \|\vrn\|_{\gamma\frac{3p}{4p-3}}^\gamma \|\vun\|_{\frac {3p}{3-p}}.
$$
Since $\frac {3p}{4p-3} = \frac{6m}{5m-1+\alpha}$, we end up with
$$
\gamma(1+\Theta) \geq \gamma \frac{6m}{5m-1+\alpha}
$$
which reduces to two conditions, using the upper bound on $1+\Theta$ presented above:
$$
\gamma \geq \frac{5m-1+\alpha}{3(m-1+\alpha)}, \qquad m \geq 1.
$$
Again, the convergence (to, however, a term with $\overline{\vr^\gamma}$) requires strict inequalities.

Finally
$$
|I_{5,2}^n| \leq C \int_\Omega \vrn\vtn|\vun|\,\dx \leq C \|\vrn\|_{\gamma(1+\Theta)} \|\vtn\|_{3m} \|\vun\|_{\frac{3p}{3-p}},
$$
provided
$$
\frac{1}{3m} + \frac{1}{\gamma(1+\Theta)} + \frac{m+1-\alpha}{6m}\leq 1.
$$  
This results into
$$
\gamma(1+\Theta) \geq \frac{6m}{5m-3+\alpha}
$$ 
and using the restrictions on $1+\Theta$ from above we end up with
$$
\gamma \geq \frac{3m+1+\alpha}{5m-3+\alpha}, \qquad \gamma \geq \frac{7m-1-\alpha}{5m-3+\alpha}.
$$

Therefore, the limit passage in the total energy balance requires
$$
\begin{aligned}
&\gamma > \frac 53, \qquad \gamma > \frac{3m+1+\alpha}{3(m+1-\alpha)}, \qquad \gamma > \frac{5m-1+\alpha}{3(m-1+\alpha)}, \\
&\gamma > \frac{3m+1+\alpha}{5m-3+\alpha}, \qquad \gamma > \frac{7m-1-\alpha}{5m-3+\alpha}, \qquad m>1.
\end{aligned}
$$
It can be easily checked that the only conditions which come into play (i.e., the most restrictive ones) are
$$
 \gamma > \frac{5m-1+\alpha}{3(m-1+\alpha)}, \qquad m>1.
$$
Furthermore,
$$
\frac{5m-1+\alpha}{3(m-1+\alpha)} > \max \Big\{\frac{3m}{2m-1+\alpha}, \frac{3m(1+\alpha)+ 2(m+1-\alpha)}{3m(1+\alpha)}\Big\}.
$$

Summarizing, to pass to the limit in the weak formulation of the continuity equation, momentum equation and entropy inequality, we need ($\Theta$ is always positive in this situation)
\begin{equation} \label{cond_var}
\gamma > \max \Big\{\frac{3m}{2m-1+\alpha}, \frac{3m(1+\alpha)+ 2(m+1-\alpha)}{3m(1+\alpha)}\Big\}, \qquad m > \max\Big\{\frac{1+\alpha}{3}, \frac{1-\alpha}{2}\Big\}.
\end{equation}
However, to pass to the limit in the total energy balance, we need (and these conditions are more restrictive than \eqref{cond_var}, thus they ensure also the limit passages in the other formulas)
\begin{equation} \label{cond_weak}
 \gamma > \frac{5m-1+\alpha}{3(m-1+\alpha)}, \qquad m>1.
\end{equation}
To conclude, we need to remove the bars over the limits of the pressure, internal energy and entropy. This is equivalent with the strong convergence of the density which we shall study in the next subsection.

\subsection{Limit passage II: strong convergence of the density}

Below, we closely follow the procedure from \cite{NP1}; some steps are identical, however, there are steps which differ due to the different form of the viscosity.

We introduce the density cut-off function
$$
T_k(z):= k T\Big(\frac{z}{k}\Big)
$$
with $T\in C^\infty([0,\infty))$ increasing such that
$$
T(z)= \left\{ \begin{array}{rl}
z & 0\leq z\leq 1 \\
\text{concave,} & 1\leq z\leq 2 \\
2 & z\geq 3
\end{array}
\right.
$$
and show that it holds
\begin{equation} \label{eff_vis_flux}
\begin{aligned}
& \overline{p(\vr,\vt) T_k(\vr)} - \Big(\frac 43 \mu(\vt) + \xi(\vt)\Big) \overline{T_k(\vr)\Div \vu} \\
= & \overline{p(\vr,\vt)} \, \overline{T_k(\vr)} - \Big(\frac 43 \mu(\vt) + \xi(\vt)\Big) \overline{T_k(\vr)}\Div \vu.
\end{aligned}
\end{equation}
This equality, called the effective viscous flux identity, can be achieved by testing the weak formulation of the momentum balance for the sequence $(\vrn,\vun,\vtn)$ by $\zeta \nabla \Delta^{-1} [1_\Omega T_k(\vrn)]$ and the limit equation from the previous section by $\zeta \nabla \Delta^{-1} [1_\Omega \overline{T_k(\vr)}]$, where $\zeta$ is an arbitrary smooth compactly supported function in $\Omega$. To show \eqref{eff_vis_flux}, we need to employ Lemmata \ref{Comm_1} and \ref{Comm_2} and perform the limit passage in lower order terms as well as suitably use several times integration by parts. We also use here the Lipschitz continuity of the shear stress tensor.

To conclude the strong convergence of the density, we also need the renormalized form of the continuity equation. Its validity is not immediate in our situation. One possibility is to use the weak formulation of the continuity equation for the limit functions and apply the Friedrich commutator lemma (Lemma \ref{Fried}). This leads to the requirement that the limit density belongs to the $L^{p'}$-space. Since $p'=\frac{p}{p-1} = \frac{6m}{3m-1+\alpha}$, we in fact need
$$
\gamma(1+\Theta) \geq \frac{6m}{3m-1+\alpha}
$$
which, using the lower bounds on $1+\Theta$, leads to restrictions
\begin{equation} \label{gamma_a}
\gamma \geq \frac{3m+1+\alpha}{3m-1+\alpha}, \qquad \gamma \geq \frac{5m+1-\alpha}{3m-1+\alpha}.
\end{equation}
This poses extra requirements in the case of variational entropy solution (and, moreover, $\frac{5m+1-\alpha}{3m-1+\alpha} \geq \frac{3m+1+\alpha}{3m-1+\alpha}$ for $m\geq \alpha$) while in the case of weak solution (i.e., it is possible to consider the weak formulation of the total energy balance) any our solution from the previous section satisfies the renormalized continuity equation since $\frac{5m-1+\alpha}{3(m-1+\alpha)}$ is larger than any of the constraint in \eqref{gamma_a}.

From now on we work only with the case of variational entropy solutions. In case the validity of the renormalized continuity equation is not ensured, we proceed in the spirit of the idea of E. Feireisl using the control of the oscillation defect measure. To this aim, we introduce (the oscillation defect measure)
\begin{equation} \label{osc_def_mea}
{\rm osc}_q[\vrn\to \vr](\Omega) = \sup_{k\geq 1} \Big(\limsup_{n\to\infty} \int_\Omega |T_k(\vrn)-T_k(\vr)|^q \,\dx \Big)
\end{equation}
with $T_k$ as above and recall that if $\vrn \rightharpoonup \vr$ in $L^1(\Omega)$, $\vun\rightharpoonup \vu$ in $W^{1,r}(\Omega)$ for some $r>1$, where the pair $(\vrn,\vun)$ is the renormalized solution to the continuity equation for any $n\in \N$, and  ${\rm osc}_q[\vrn\to \vr](\Omega) <\infty$ for $\frac 1q < 1-\frac 1r$, then the limit pair $(\vr,\vu)$ also satisfies the renormalized continuity equation (see, e.g., \cite[Lemma 3.8]{FNbook}). Recall now \cite[Lemma 18]{NP1} which says that (the positive number $d$ comes either from \eqref{pre_prop2} or \eqref{pre_four})
$$
\begin{aligned} 
d\limsup_{n\to \infty} \int_\Omega |T_k(\vrn)-T_k(\vr)|^{\gamma+1} \,\dx &\leq \int_\Omega (\overline{p(\vr,\vt)T_k(\vr)} - \overline{p(\vr,\vt)} \, \overline{T_k(\vr)})\, \dx \\
d\limsup_{n\to \infty} \int_\Omega \frac{1}{(1+\vt)^\alpha}|T_k(\vrn)-T_k(\vr)|^{\gamma+1} \,\dx &\leq \int_\Omega (\overline{p(\vr,\vt)T_k(\vr)} - \overline{p(\vr,\vt)} \, \overline{T_k(\vr)})\, \dx.
\end{aligned}
$$
The proof is performed for $\alpha=1$, however, the modification for $0\leq \alpha <1$ is straightforward here. Therefore we have
$$
\begin{aligned}
d \limsup_{n\to \infty}& \int_\Omega |T_k(\vrn)-T_k(\vr)|^{\gamma+1} \,\dx \leq d\limsup_{n\to \infty} \int_\Omega (\vrn^\gamma-\vr^\gamma)(T_k(\vrn)-T_k(\vr)) \,\dx \\
&= d \int_\Omega \Big(\overline{\vr^\gamma T_k(\vr)} -\overline{\vr^\gamma} \,\overline{T_k(\vr)}\Big)\,\dx + d \int_\Omega (\vr^\gamma -\overline{\vr^\gamma}) ( T_k(\vr) -\overline{T_k(\vr)})\,\dx.
\end{aligned}
$$
Setting $G(x,z) = d |T_k(z)-T_k(\vr(x))|^{\gamma+1}$ we have (we employ here Lemma \ref{weak_convex} to show that the last term on the right hand-side is non-positive)
$$
\overline{G(\cdot,\vr)} \leq \overline{p(\vr,\vt)} \, \overline{T_k(\vr)} - \overline{p(\vr,\vt)T_k(\vr)} 
$$
and due to the effective viscous flux identity \eqref{eff_vis_flux}
$$
\overline{G(\cdot,\vr)} \leq \Big(\frac{4}{3}\mu(\vt) +\xi(\vt)\Big) (\overline{T_k(\vr) \Div \vu} - \overline{T_k(\vr)}\Div\vu).
$$
Then
$$
\begin{aligned}
\int_\Omega(1+\vt)^{-\alpha} \overline{G(\cdot,\vr)} \,\dx &\leq C \sup_{n\in \N} \|\Div \vun\|_p \limsup_{n\to \infty}\|T_k(\vrn)-T_k(\vr)\|_{p'} \\
&\leq C \limsup_{n\to \infty}\|T_k(\vrn)-T_k(\vr)\|_{p'}
\end{aligned}
$$
with $p' = \frac{6m}{3m-1+\alpha}$. On the other hand
$$
\begin{aligned}
\int_\Omega &|T_k(\vrn)-T_k(\vr)|^q\,\dx \leq \int_\Omega |T_k(\vrn)-T_k(\vr)|^q (1+\vt)^{-\frac{q\alpha}{\gamma+1}} (1+\vt)^{\frac{q\alpha}{\gamma+1}}\,\dx \\
&\leq \Big( \int_\Omega (1+\vt)^{-\alpha} |T_k(\vrn)-T_k(\vr)|^{\gamma+1}\,\dx \Big)^{\frac{q}{\gamma+1}} \Big(\int_\Omega (1+\vt)^{\frac{q\alpha}{\gamma+1-q}}\,\dx \Big)^{\frac{\gamma+1-q}{\gamma+1}}.
\end{aligned}
$$
To conclude, we need first $\frac{q\alpha}{\gamma+1-q}\leq 3m$, i.e. $q\leq \frac{3m(\gamma+1)}{3m+\alpha}$ and further $q>p'=\frac{6m}{3m-1+\alpha}$ which leads to
$$
\frac{3m(\gamma+1)}{3m+\alpha} > \frac{6m}{3m-1+\alpha}.
$$
This results into
$$
\gamma > \frac{3m+1+\alpha}{3m-1+\alpha}.
$$
Whence the renormalized continuity equation holds for the variational entropy solution for $\gamma > \frac{3m+1+\alpha}{3m-1+\alpha}$. Note that this condition is precisely the same as one of the conditions in \eqref{gamma_a}. Therefore we may relax the condition $\gamma \geq \frac{5m+1-\alpha}{3m-1+\alpha}$ on the validity of the renormalized continuity equation by Friedrich commutator lemma and consider just the condition coming from the estimate of the oscillation defect measure.

The strong convergence of the density then follows directly, exactly as in \cite[pp. 312--313]{NP1} in case of the oscillation defect measure estimate or  \cite[pp. 12--13]{P} in case when the renormalized continuity equation holds directly. Theorem \ref{t_heat_flux} is proved.

\section{Proof for the Dirichlet boundary conditions for the temperature}
\label{5}
 
In this section we have to restrict ourselves to the pressure form \eqref{pre_two} together with the assumption that the Third law of thermodynamics holds. In particular, the entropy satisfies bounds \eqref{ent_four}.

\subsection{Energy and entropy estimates}
\label{5.1}

Let $(\vrn,\vun,\vtn)$ be a sequence of weak solutions to our problem \eqref{NSF_1}, \eqref{vel_Dir} (or \eqref{vel_Nav}), \eqref{temp_Dir} and \eqref{tot_mass} in the sense of Definition \ref{d2}, formed, e.g., by a convergent sequence of total masses $M_n \to M>0$. We proceed similarly as in Section \ref{4}, however, this time we have at our disposal only the ballistic energy inequality \eqref{ball_ene} instead of the entropy one. In particular, the power of densities on the right hand-side of the estimates will be higher than in the previous section which is caused by the presence of density, velocity and temperature on the right hand-side of the ballistic energy inequality (recall that on the right hand-side of the entropy inequality, there was just a constant). 

Recall that we have $\psi_D$ the harmonic extension of the boundary data for the temperature
$$
\begin{aligned}
\Delta \psi_D &= 0 \qquad \text{ in } \Omega \\
\psi_D &= \vt_D \qquad \text { at } \partial \Omega.
\end{aligned}
$$
Recall also that due to our assumptions $\psi_D \in W^{2,q}(\Omega) \hookrightarrow C^1(\overline{\Omega})$ for $q>3$. Note that if $\Omega \in C^{1,1}$ and $\vt_D \in W^{2-\frac 1q,q}(\partial \Omega)$, then the smoothness of $\psi_D$ follows. Moreover, by maximum principle, if $0 <C_1 \leq \vt_D \leq C_2$ (for the continuous representative on $\partial \Omega$), then $C_1 \leq \psi_D \leq C_2$ (again for the continuous representative now in $\Omega$).

The ballistic energy inequality reads as follows (for the Navier boundary conditions we take for simplicity $\alpha =0$, i.e., we consider the total slip boundary conditions)
\begin{equation}\label{ball_ene_n}
\begin{aligned}
\int_\Omega &\Big(\frac{\tn{S}(\vtn,\nabla \vun):\nabla \vun}{\vtn} + \frac{\kappa(\vtn)|\nabla \vtn|^2}{\vt^2}\Big)\psi_D \dx   \\
\leq &\int_\Omega \Big(\vrn \vc{f} \cdot \vun - \vrn s(\vrn\vtn)\vun \cdot \nabla \psi_D + \frac{\kappa(\vtn) \nabla \vtn \cdot \nabla \psi_D}{\vtn}\Big)\dx. 
\end{aligned}
\end{equation}

We proceed as in \cite{P} (see also \cite{FC}). We have for $K(z):= \int_1^x \frac{\kappa(s)}{s} \, {\rm d}s$
$$
\begin{aligned}
\int_\Omega \frac{\kappa(\vtn) \nabla \vtn \cdot \nabla \psi_D}{\vtn} \dx &= \int_\Omega \nabla K(\vtn)\cdot \nabla \psi_D \dS \\
& = -\int_\Omega K(\vtn) \Delta \psi_D \dx 
+ \int_{\partial \Omega} K(\vtn) \pder{\psi_D}{\vc{n}}\dS \\
& = \int_{\partial \Omega} K(\psi_D) \pder{\psi_D}{\vc{n}}\dS
\end{aligned}
$$
and the term is bounded independently of $n$ by the data of the problem. Next
$$
\Big|\int_\Omega \vrn\vc{f}\cdot \vun \dx\Big| \leq C \|\vun\|_{\frac{3p}{3-p}} \|\vrn\|_{\frac{3p}{4p-3}}
$$
and
$$
\begin{aligned}
\Big|\int_\Omega \vrn s(\vrn,\vtn)\vun \cdot \nabla \psi_D\dx\Big| &\leq C \int_\Omega \big(\vrn + \vrn [\log\vrn]^+ + \vrn[\log \vtn]^+\big)|\vun| \dx \\
&\leq C(q,\delta) \|\vun\|_{\frac{3p}{3-p}} \big(\|\vrn\|_q^{1+\delta} +\|\vrn\|_q (1+ \|[\log \vtn]^+\|_r)\big), 
\end{aligned}
$$
where $q> \frac{3p}{4p-3}$, but can be taken arbitrarily close to it, $\delta >0$, can be  taken arbitrarily small and
$$
\frac 1q + \frac{3-p}{3p} + \frac 1r =1.
$$
Thus, using \eqref{48a} we get for $p = \frac{6m}{3m+1-\alpha}$
$$
\begin{aligned}
\|\vun\|_{1,p}^p &\leq C \Big(\int_\Omega \frac{{\tn{S}(\vtn,\nabla \vun):\nabla \vun}}{\vtn}\dx\Big)^{\frac p2} \Big(\int_\Omega \vtn ^{\frac{1-\alpha}{2}\, \frac{2p}{2-p}}\dx\Big)^{\frac {2-p}{2}} \\
&\leq
C \|\vun\|_{1,p}^{\frac p2}\|\vtn\|_{3m}^{\frac{1-\alpha}{2}p} \big(1+\|\vrn\|_q^{1+\delta} + \|\vrn\|_q (1+\|[\log \vtn]^+\|_r)\big)^{\frac p2}
\end{aligned}
$$
which leads to
$$
\|\vun\|_{1,p} \leq C \|\vtn\|_{3m}^{1-\alpha} \big(1+\|\vrn\|_q^{1+\delta} + \|\vrn\|_q (1+\|\vtn\|_{3m}^\delta)\big)
$$
as well as to
$$
\begin{aligned}
\|\vtn\|_{3m}^m &\leq C\Big(1+ \int_\Omega \big|\nabla \big(\vtn^{\frac m2}\big)\big|^2 \Big)\dx \\
& \leq C (1+\|\vun\|_{1,p})  \big(1+ \|\vrn\|_q^{1+\delta} + \|\vrn\|_q (1+\|[\log \vtn]^+\|_r)\big) \\
&\leq C \big(1+ \|\vtn\|_{3m}^{1-\alpha +\delta}\big) \big(1+\|\vrn\|_q^{1+\delta}\big)^2,
\end{aligned}
$$ 
where the role $q$ and $\delta$ is the same as above.
Thus, provided $m>1-\alpha+\delta$, 
$$
\begin{aligned}
\|\vtn\|_{3m} & \leq C\Big(1+ \|\vrn\|_{q}^{\frac{2(1+\delta)}{m-1+\alpha-\delta}}\Big)\\
\|\vun\|_{1,p} &\leq C \Big(1+\|\vrn\|_q^{\frac{(m+1-\alpha-\delta)(1+\delta)}{m-1-\delta +\alpha}}\Big).
\end{aligned}
$$
Since $q> \frac{3p}{4p-3}$ and can be chosen arbitrarily close to this number and $\delta >0$ can be chosen arbitrarily close to $0$, we set in the following computations formally $q= \frac{3p}{4p-3}$ and $\delta =0$. After many computations we obtain at the end several inequalities and we replace those
which include equality by the strict ones. This will finally provide conditions under which the estimates are available and finally also our solution exists. 

Whence, we consider
\begin{equation} \label{est_5_1}
\begin{aligned}
\|\vtn\|_{3m} & \leq C\Big(1+ \|\vrn\|_{\frac{3p}{4p-3}}^{\frac{2}{m-1+\alpha}}\Big)\\
\|\vun\|_{1,p} &\leq C \Big(1+\|\vrn\|_{\frac{3p}{4p-3}}^{\frac{m+1-\alpha}{m-1+\alpha}}\Big),
\end{aligned}
\end{equation}
provided $m>1-\alpha$. Furthermore, if $\frac{3p}{4p-3} = \frac{6m}{5m-1+\alpha} < \gamma(1+\Theta)$, we have (recall that $\frac{6m}{5m-1+\alpha} >1$ for $m>0$ and $\alpha \in [0,1]$)
\begin{equation} \label{int_5}
\|\vrn\|_{\frac{3p}{4p-3}} \leq C(M) \|\vrn\|_{\gamma(1+\Theta)}^{\frac{\gamma(1+\Theta)(m+1-\alpha)}{6m(\gamma(1+\Theta)-1)}}.
\end{equation}

\subsection{Density estimates by Bogovskii operator}\label{s:52}

We now proceed exactly as in the previous section for the heat flux boundary conditions. We use as a test function in the weak formulation of the momentum equation the function
$$
\vcg{\varphi} := \mathcal B \Big(p(\vrn,\vtn)^\Theta - \frac{1}{|\Omega|} \int_\Omega p(\vrn,\vtn)^\Theta \dx \Big)
$$
where $\mathcal B$ is the Bogovskii operator (see Lemma \ref{Bogovskii}), $\Theta >0$. 

After exactly the same calculations as in the previous section we show that the right-hand side of \eqref{Bog_test} is finite provided (recall that all equalities are changed into strict inequalities due to small factors in the estimate)
\begin{equation} \label{est_theta5}
\Theta < \min\Big\{\frac{3m-1-\alpha}{3m+1+\alpha}, \frac{\gamma(2m-1+\alpha)-3m}{\gamma(m+1-\alpha)}\Big\}
\end{equation}
and
\begin{equation} \label{est_theta5a}
m > \frac{1+\alpha}{3}, \qquad \gamma > \frac{3m}{2m-1+\alpha}.
\end{equation}
Moreover, exactly as in the previous section it follows that the condition needed in the interpolation inequality for the density
$$
\gamma(1+\Theta) > \frac{6m}{5m+1-\alpha}
$$
can be fulfilled provided $\Theta$ satisfies \eqref{est_theta5}. However, the upper bounds on $\Theta$ under which the application of the Bogovskii operator above really provides estimate of the density differs from Section 4. We have
$$
\begin{aligned}
|I_1| & \leq \int_\Omega \vrn |\vun|^2 |\nabla \vcg{\varphi}|\dx \leq \varepsilon \|p(\vrn,\vtn)\|_{1+\Theta}^{1+\Theta} + C(\varepsilon) \|\vrn\|_{\gamma(1+\Theta)}^{1+\Theta} \|\vun\|_{\frac{3p}{3-p}}^{2(1+\Theta)} \\
&\leq  \varepsilon \|p(\vrn,\vtn)\|_{1+\Theta}^{1+\Theta} + C(\varepsilon) \|\vrn\|_{\gamma(1+\Theta)}^{1+\Theta + 2(1+\Theta) \frac{\gamma(1+\Theta)(m+1-\alpha)}{6m(\gamma(1+\Theta)-1)}\frac{(m+1-\alpha)}{m-1+\alpha}}.
\end{aligned}
$$
Above, we used \eqref{est_5_1} as well as \eqref{int_5}. Whence we require
$$
\gamma(1+\Theta) > 1+\Theta + 2(1+\Theta) \frac{\gamma(1+\Theta)(m+1-\alpha)}{6m(\gamma(1+\Theta)-1)}\frac{(m+1-\alpha)}{m-1+\alpha}
$$
which after standard computations results into
\begin{equation} \label{Theta_lower1}
\Theta > \frac{\gamma(m+1-\alpha)^2 -3m(\gamma-1)^2 (m-1+\alpha)}{\gamma\big(3m(\gamma-1)(m-1+\alpha) -(m+1-\alpha)^2\big)}.
\end{equation}
Below, we shall show that $\gamma >\frac{5m-1+\alpha}{3(m-1+\alpha)}$ and under this assumption the denominator is always positive provided $m>1-\alpha$.

Next we have
$$
\begin{aligned}
|I_2| &\leq \int_\Omega \big|\tn{S}(\vtn,\nabla \vun):\nabla \vcg{\varphi}\big|\dx \leq \big(1+\|\vtn\|_{3m}^\alpha\big) \|\nabla \vun\|_p \|p(\vrn,\vtn)\|_{1+\Theta}^\Theta \\
& \leq \varepsilon \|p(\vrn,\vtn)\|_{1+\Theta}^{1+\Theta}  + C(\varepsilon) (1+\|\vtn\|_{3m}^{\alpha(1+\Theta)}) \|\nabla \vun\|_p^{1+\Theta} \\
& \leq \varepsilon \|p(\vrn,\vtn)\|_{1+\Theta}^{1+\Theta}  + C(\varepsilon) \Big(1+ \|\vrn\|_{\gamma(1+\Theta)}^{\frac{2\alpha(1+\Theta)}{m-1+\alpha}\frac{\gamma(1+\Theta)(m+1-\alpha)}{6m(\gamma(1+\Theta) -1)} + \frac{(1+\Theta)(m+1-\alpha) \gamma(1+\Theta)(m+1-\alpha)}{(m-1+\alpha) 6m(\gamma(1+\Theta) -1)}}\Big).
\end{aligned}
$$
Therefore we require
$$
\gamma (1+\Theta) > (1+\Theta) \frac{\gamma(1+\Theta)(m+1-\alpha) (m+1+\alpha)}{(m-1+\alpha) 6m(\gamma(1+\Theta) -1)}
$$
which leads to
\begin{equation} \label{Theta_lower_2}
\Theta > \frac{(m+1)^2 -\alpha^2 -6m(\gamma-1)(m-1+\alpha)}{6m\gamma (m-1+\alpha) -(m+1)^2+\alpha^2}.
\end{equation}
As above, for $\gamma >\frac{5m-1+\alpha}{3(m-1+\alpha)}$ the denominator is always positive provided $m>\max\{1-\alpha,\frac{1+\alpha}{3}\}$.
Finally we consider (the term $\int_\Omega \vrn\vc{f}_n\cdot \vcg{\varphi}_n\dx$ is of lower order)
$$
I_4 = \int_\Omega p(\vrn,\vtn)\dx \frac{1}{|\Omega|}\int_\Omega p(\vrn,\vtn)^\Theta \dx.
$$
Recall that similarly as in the previous section the most restrictive term is
$$
I_{4,4} = C\int_\Omega \vrn\vtn\dx \int_\Omega (\vrn\vtn)^\Theta \dx,
$$
hence as in the previous section
$$
\begin{aligned}
|I_{4,4}| &\leq \varepsilon \|p(\vrn,\vtn)\|_{1+\Theta}^{1+\Theta} + C(M,\varepsilon) \|\vtn\|_{3m}^{1+\Theta} \\
&\leq \varepsilon \|p(\vrn,\vtn)\|_{1+\Theta}^{1+\Theta} + C(M,\varepsilon) \|\vrn\|_{\gamma(1+\Theta)}^{\frac{2(1+\Theta)}{m-1+\alpha}\frac{\gamma(1+\Theta)(m+1-\alpha)}{6m(\gamma(1+\Theta)-1)}}.
\end{aligned}
$$
Therefore
$$
\gamma (1+\Theta) > \frac{(1+\Theta)}{m-1+\alpha}\frac{\gamma(1+\Theta)(m+1-\alpha)}{3m(\gamma(1+\Theta)-1)}
$$
which leads to
$$
1> \frac{(1+\Theta)}{m-1+\alpha}\frac{(m+1-\alpha)}{3m(\gamma(1+\Theta)-1)}.
$$
However, this bound is clearly less restrictive than the bound coming from the estimate of $I_2$ as
$$
\frac{(1+\Theta)(m+1-\alpha) (m+1+\alpha)}{(m-1+\alpha) 6m(\gamma(1+\Theta) -1)} > \frac{(1+\Theta)}{m-1+\alpha}\, \frac{(m+1-\alpha)(\gamma(1+\Theta)-1)}{3m (\gamma(1+\Theta)-1)};
$$
this inequality leads after easy computations to $m+1+\alpha>2$ which is easily satisfied as $m>1-\alpha$.
Concluding this part, the estimate
\begin{equation} \label{est_final_1}
\|\vun\|_{1,p} + \|\vtn\|_{3m} + \|\nabla \vtn^{\frac m2}\|_2 + \|\vrn\|_{\gamma(1+\Theta)} \leq C
\end{equation}
for some
$$
\Theta < \min\Big\{\frac{3m-1-\alpha}{3m+1+\alpha}, \frac{\gamma(2m-1+\alpha)-3m}{\gamma(m+1-\alpha)}\Big\}
$$
provided
$$
m > \max\Big\{1-\alpha,\frac{1+\alpha}{3}\Big\}, \qquad \gamma > \frac{3m}{2m-1+\alpha}
$$
and
$$
\begin{aligned}
& \min\Big\{\frac{3m-1-\alpha}{3m+1+\alpha}, \frac{\gamma(2m-1+\alpha)-3m}{\gamma(m+1-\alpha)}\Big\} >\\
& \max \Big\{\frac{\gamma(m+1-\alpha)^2 -3m(\gamma-1)^2 (m-1+\alpha)}{\gamma\big(3m(\gamma-1)(m-1+\alpha) -(m+1-\alpha)^2\big)},\frac{(m+1)^2 -\alpha^2 -6m(\gamma-1)(m-1+\alpha)}{6m\gamma (m-1+\alpha) -(m+1)^2+\alpha^2}\Big\}.
\end{aligned}
$$
Recall also that we used $\gamma > \frac{5m-1+\alpha}{3(m-1+\alpha)}$ to show that, the denominators of the fractions on the right hand side of the above inequality are positive. However, this is at this point just a sufficient, not necessary condition. To conclude we need to verify that the last inequality does not define empty set and if this inequality requires further restrictions, to compute them.

First, let us make the following observation. The inequality
$$
\frac{\gamma(m+1-\alpha)^2 -3m(\gamma-1)^2 (m-1+\alpha)}{\gamma\big(3m(\gamma-1)(m-1+\alpha) -(m+1-\alpha)^2\big)} > \frac{(m+1)^2 -\alpha^2 -6m(\gamma-1)(m-1+\alpha)}{6m\gamma (m-1+\alpha) -(m+1)^2+\alpha^2}
$$
holds provided
$$
\gamma(m+1-3\alpha) + m+1+\alpha > 0.
$$
For $\alpha \in [0,\frac 12]$ it can be verified that the inequality holds true; if $\alpha >\frac 12$ (i.e., it includes $\alpha=1$), it may happen that for $\gamma$ large enough the inequality is not fulfilled. 

First we consider
$$
\frac{\gamma(2m-1+\alpha)-3m}{\gamma(m+1-\alpha)} > \frac{\gamma(m+1-\alpha)^2 -3m(\gamma-1)^2 (m-1+\alpha)}{\gamma\big(3m(\gamma-1)(m-1+\alpha) -(m+1-\alpha)^2\big)}.
$$
This leads, after tedious, however straightforward computation to the quadratic inequality
$$
3(m-1+\alpha)\gamma^2 -\gamma (8m-4(1-\alpha)) + 5m-1+\alpha >0.
$$
The value $\gamma =1$ is one root of the quadratic polynomial on the left-hand side and then it is not difficult to see that the other root is $\gamma =  \frac{5m-1+\alpha}{3(m-1+\alpha)}$, i.e. we get precisely the above mentioned restriction
$$
\gamma >\frac{5m-1+\alpha}{3(m-1+\alpha)},
$$
which now becomes necessary condition for the estimates to hold.

Next, the inequality
$$
\frac{3m-1-\alpha}{3m+1+\alpha} > \frac{\gamma(m+1-\alpha)^2 -3m(\gamma-1)^2 (m-1+\alpha)}{\gamma\big(3m(\gamma-1)(m-1+\alpha) -(m+1-\alpha)^2\big)}
$$
leads to quadratic inequality
\begin{equation}\label{QI1}
\begin{aligned} 
&6m\gamma^2 (m-1+\alpha) -\gamma[ 2(m+1-\alpha)^2 + (m-1+\alpha)(9m+1+\alpha)] \\
&+ (m-1+\alpha)(3m+1+\alpha) >0
\end{aligned}
\end{equation}
which holds true provided $\gamma$ is larger than $\gamma_2$, which is the larger of two roots $\gamma_{1,2}$ of the corresponding quadratic equation. These roots can be expressed as
$$
\gamma_{1,2} = \frac{11m^2 + m(6\alpha-4)+3\alpha^2-4\alpha+1 \pm \sqrt{D_1(m,\alpha)}}{12m(m+\alpha-1)}
$$
with $D_1(m,\alpha)$ being a polynomial of fourth order in $m$ and $\alpha$. If we try to determine, whether $\gamma > \gamma_2$ is more restrictive condition than $\gamma > \frac{5m-1+\alpha}{3(m-1+\alpha)}$, i.e. study the inequality $\gamma_2 < \frac{5m-1+\alpha}{3(m-1+\alpha)}$, we recover after further tedious computation the following simple inequality
$$
32m(m-1)(m-\alpha+1)(m+\alpha-1) > 0.
$$
Recalling that $m > 1-\alpha$, we conclude that for $m > 1$ the condition $\gamma > \gamma_2$ is less restrictive, while for $m \leq 1$ the condition $\gamma > \gamma_2$ is more restrictive than $\gamma > \frac{5m-1+\alpha}{3(m-1+\alpha)}$.

We may therefore conclude for $\alpha \in [0,\frac 12]$. Since
$$
\frac{5m-1+\alpha}{3(m-1+\alpha)} > \frac{3m}{2m-1+\alpha}
$$
for any $\alpha \in [0,1]$, the bounds for the velocity, density and velocity \eqref{est_final_1} holds true with
$$
\Theta < \min\Big\{\frac{3m-1-\alpha}{3m+1+\alpha}, \frac{\gamma(2m-1+\alpha)-3m}{\gamma(m+1-\alpha)}\Big\},
$$
arbitrarily close to this number, provided (recall that $1-\alpha\geq \frac{1+\alpha}{3}$ for $\alpha \in [0,\frac 12]$)
$$
\begin{aligned}
m&>1-\alpha, \qquad \gamma >  \frac{5m-1+\alpha}{3(m-1+\alpha)} \\
&6m\gamma^2 (m-1+\alpha) -\gamma[ 2(m+1-\alpha)^2 +  (m-1+\alpha)(9m+1+\alpha)] \\
&+ (m-1+\alpha)(3m+1+\alpha) >0,
\end{aligned}
$$
where we recall that the quadratic inequality represents a restriction only in the case $m < 1$.

Now, let $\frac 12 < \alpha \leq 1$. We have to study the other two inequalities. First, consider
$$
\frac{3m-1-\alpha}{3m+1+\alpha} > \frac{(m+1)^2 -\alpha^2 -6m(\gamma-1)(m-1+\alpha)}{6m\gamma (m-1+\alpha) -(m+1)^2+\alpha^2}.
$$
This reduces to
$$
\gamma > \frac{2(m+\alpha)}{3(m-1+\alpha)},
$$
however, for $m >\frac{1+\alpha}{3}$ we have
$$
\frac{5m-1+\alpha}{3(m-1+\alpha)} > \frac{2(m+\alpha)}{3(m-1+\alpha)}
$$
and thus the inequality above does not lead to any further restrictions.

Last we consider 
$$
\frac{\gamma(2m-1+\alpha)-3m}{\gamma(m+1-\alpha)} > \frac{(m+1)^2 -\alpha^2 -6m(\gamma-1)(m-1+\alpha)}{6m\gamma (m-1+\alpha) -(m+1)^2+\alpha^2}.
$$
This can be easily transformed into the quadratic inequality
$$
6m\gamma^2 (m-1+\alpha) -\gamma[2(m-1+\alpha)(4m+1-\alpha) + (m+1)^2 -\alpha^2] + (m+1)^2-\alpha^2>0.
$$
We claim that this inequality does not represent any further restriction. Indeed, it yields similarly as above a condition
$$
\gamma > \gamma_4 = \frac{9m^2+m(6\alpha-4)-3\alpha^2+4\alpha-1 + \sqrt{D_2(m,\alpha)}}{12m(m+\alpha-1)}
$$
with $D_2(m,\alpha)$ being a polynomial of a fourth order in $m$ and $\alpha$. The inequality $\gamma_4 < \frac{5m-1+\alpha}{3(m-1+\alpha)}$ then simplifies into
$$
32m(2m-1)(m-\alpha+1)^2 > 0,
$$
which is always satisfied since for $\alpha > \frac12$ we have $m > \frac{1+\alpha}{3} > \frac12$.

Concluding, estimate \eqref{est_final_1} holds in case of $\frac 12 <\alpha \leq 1$ provided
$$
\Theta < \min\Big\{\frac{3m-1-\alpha}{3m+1+\alpha}, \frac{\gamma(2m-1+\alpha)-3m}{\gamma(m+1-\alpha)}\Big\},
$$ 
arbitrarily close to this number,
and 
$$
\begin{aligned}
& m > \frac{1+\alpha}{3}, \qquad \gamma > \frac{5m-1+\alpha}{3(m-1+\alpha)}, \\
&6m\gamma^2 (m-1+\alpha) -\gamma[ 2(m+1-\alpha)^2 +  (m-1+\alpha)(9m+1+\alpha)] \\
&+ (m-1+\alpha)(3m+1+\alpha) >0, 
\end{aligned}
$$
where we again recall that the quadratic inequality represents a restriction only for $m < 1$.

\subsection{Limit passage for the sequence of solutions}

To pass to the limit $n \to \infty$ in our problem, we proceed exactly as in the previous section. Recall that we need to have the renormalized continuity equation which can be shown to hold (see Subsection 4.4) provided
$$
\gamma > \frac{3m+1+\alpha}{3m-1+\alpha}.
$$
If we consider the inequality
$$
\frac{5m-1+\alpha}{3(m-1+\alpha)} > \frac{3m+1+\alpha}{3m-1+\alpha},
$$
we easily check that it is equivalent to
\begin{equation} \label{q_ineq_1}
6m^2-2m-4m\alpha +(1-\alpha)^2 + 3(1-\alpha^2) >0.
\end{equation}
If $0\leq \alpha \leq \frac12$,  it is not difficult to see that the left-hand side is monotone in $m$ for $m>1-\alpha$ and thus, plugging in $m=1-\alpha$ we conclude that inequality \eqref{q_ineq_1} holds in this range of $\alpha$ always, provided $m >1-\alpha$. Similar argument applies also in the case $\frac 12 <\alpha\leq 1$ and $m > 1$, where \eqref{q_ineq_1} also holds. Therefore, in both these cases $\gamma > \frac{3m+1+\alpha}{3m-1+\alpha}$ does not represent a further restriction. This however does not apply for $\frac 12 <\alpha\leq 1$ and $m<1$. 

Note that for $m < 1$ the more restrictive condition than $\gamma > \frac{5m-1+\alpha}{3(m-1+\alpha)}$ is the quadratic inequality \eqref{QI1}. However, comparing the bound arising from \eqref{QI1} with $\frac{3m+1+\alpha}{3m-1+\alpha}$ does not yield a simple inequality as was the case in previous calculations in Section \ref{s:52}. Instead, one obtains
$$
3m^3 + m^2(2-5\alpha)+m(\alpha^2-4\alpha+3)+\alpha(1-\alpha)^2 < 0.
$$
Plugging in $m = \frac{1+\alpha}{3}$ and replacing inequality sign with equality sign we obtain cubic equation
$$
2\alpha^3-8\alpha^2+2\alpha+3 = 0.
$$
This equation has three real roots, one of them (denoted by $\alpha_0$) belonging to the interval $(0,1)$. Its value is approximately $\alpha_0 \approx 0.87$, however it does not have a simple algebraic expression. We can state that for $\alpha \in (\alpha_0,1]$ there exists $1 > m_\alpha > \frac{1+\alpha}{3}$ such that for $m \in (\frac{1+\alpha}{3},m_\alpha)$ the condition $\gamma > \frac{3m+1+\alpha}{3m-1+\alpha}$ is more restrictive than the quadratic inequality \eqref{QI1}.

Concluding, the variational entropy ballistic energy solution to our problem exists, provided 
\begin{itemize}
\item[a)] $\alpha \in [0,1]$, $m > 1$ and $\gamma >  \frac{5m-1+\alpha}{3(m-1+\alpha)}$,
\item[b)] $\alpha \in [0,\frac12]$, $m \in (1-\alpha,1]$ and 
\begin{equation}\label{QI1n}
\begin{aligned} 
&6m\gamma^2 (m-1+\alpha) -\gamma[ 2(m+1-\alpha)^2 + (m-1+\alpha)(9m+1+\alpha)] \\
&+ (m-1+\alpha)(3m+1+\alpha) >0,
\end{aligned}
\end{equation}
\item[c)] $\alpha \in (\frac12,1]$, $m \in (\frac{1+\alpha}{3},1]$ and both \eqref{QI1n} and $\gamma > \frac{3m+1+\alpha}{3m-1+\alpha}$ are satisfied.
\end{itemize}


Finally we need to verify under which conditions also the weak formulation of the total energy balance holds. We in fact only need to check under which conditions all the terms are bounded and then change the possible equalities into strict inequalities. Since the computations are exactly the same as in the case of heat flux boundary conditions, we only recall the results from the previous section. Therein, it was shown that we have to require that
$$
m>1, \qquad \gamma > \frac{5m-1+\alpha}{3(m-1+\alpha)}.
$$
This simply eliminates possibilities b) and c) above.
Theorem \ref{t_Dirichlet_bc} is proved.

\section*{Acknowledgment}
The work of Ond\v rej Kreml was supported by the Institute of Mathematics, Czech Academy of Sciences (RVO 67985840). Ond\v rej Kreml, Milan Pokorn\'y and Emil Sk\v{r}\'\i\v{s}ovsk\'y were supported by the Czech Science Foundation, project No. GA22-01591S. The work of Tomasz Piasecki was supported by the National Science Centre (NCN) project 2022/45/B/ST1/03432.

\end{document}